\documentclass[11pt]{article}
\usepackage{amssymb}
\usepackage[usenames,dvipsnames]{xcolor}
\usepackage{amsthm}
\usepackage{amsmath}
\usepackage{amsfonts}

\usepackage{graphicx}

 \theoremstyle{definition}

 \numberwithin{equation}{section}

\newtheorem{theorem}{Theorem}[section]

\theoremstyle{definition}

\theoremstyle{remark}
\newtheorem{remark}[theorem]{Remark}

\begin{document}
\title{A new formulation for the 3-D Euler equations with an application to subsonic flows in a cylinder}
\author{Shangkun WENG\footnote{Pohang Mathematics Institute,
Pohang University of Science and Technology. Hyoja-Dong San 31, Nam-Gu
Pohang, Gyungbuk 790-784, Korea. {\it Email: skweng@postech.ac.kr, skwengmath@gmail.com}}}
\date{}
\maketitle

\begin{abstract}

In this paper, a new formulation for the three dimensional Euler
equations is derived. Since the Euler system is
hyperbolic-elliptic coupled in a subsonic region, so an effective decoupling of the hyperbolic and elliptic modes is essential for any development of the theory. The key idea in our formulation is to use the Bernoulli's law to reduce the dimension of the velocity field by defining new variables $(1,\beta_2=\frac{u_2}{u_1},\beta_3=\frac{u_3}{u_1})$ and replacing $u_1$ by the Bernoulli's function $B$ through $u_1^2=\frac{2(B-h(\rho))}{1+\beta_2^2+\beta_3^2}$. We find a conserved quantity for flows with a constant Bernoulli's function, which behaves like the scaled vorticity in the 2-D case. More surprisingly, a system of new conservation laws can be derived, which is new even in the two dimensional case. We use this new formulation to construct a smooth subsonic Euler flow in a rectangular cylinder, which is also required to be adjacent to some special subsonic states. The same idea can be applied to obtain similar information for the 3-D incompressible Euler equations, the self-similar Euler equations, the steady Euler equations with damping, the steady Euler-Poisson equations and the steady Euler-Maxwell equations.

\end{abstract}

\begin{center}
\begin{minipage}{5.5in}
2010 Mathematics Subject Classification:35Q31; 35Q35;76G25.

\

Key words: Subsonic flow, Bernoulli's law, Conservation laws,
Hyperbolic-elliptic coupled.
\end{minipage}
\end{center}

\

\everymath{\displaystyle}
\newcommand {\eqdef }{\ensuremath {\stackrel {\mathrm {\Delta}}{=}}}

% math symbols
\def\Xint #1{\mathchoice
{\XXint \displaystyle \textstyle {#1}} %
{\XXint \textstyle \scriptstyle {#1}} %
{\XXint \scriptstyle \scriptscriptstyle {#1}} %
{\XXint \scriptscriptstyle \scriptscriptstyle {#1}} %
\!\int}
\def\XXint #1#2#3{{\setbox 0=\hbox {$#1{#2#3}{\int }$}
\vcenter {\hbox {$#2#3$}}\kern -.5\wd 0}}
\def\ddashint {\Xint =}
\def\dashint {\Xint -}
\def\clockint {\Xint \circlearrowright } % GOOD !
\def\counterint {\Xint \rotcirclearrowleft } % Good for Computer Modern !
\def\rotcirclearrowleft {\mathpalette {\RotLSymbol { -30}}\circlearrowleft }
\def\RotLSymbol #1#2#3{\rotatebox [ origin =c ]{#1}{$#2#3$}}

\def\aint{\dashint}

\def\arraystretch{2}
\def\eps{\varepsilon}

\def\s#1{\mathbb{#1}} % set
\def\t#1{\tilde{#1}} %new variables
\def\b#1{\overline{#1}}
\def\N{\mathcal{N}} %Nozzle
\def\M{\mathcal{M}} %Mach number
\def\R{{\mathbb{R}}}

\def\ba{\begin{array}}
\def\ea{\end{array}}
\def\be{\begin{equation}}
\def\ee{\end{equation}}

\def\bes{\begin{mysubequations}}
\def\ees{\end{mysubequations}}

\def\cz#1{\|#1\|_{C^{0,\alpha}}}
\def\ca#1{\|#1\|_{C^{1,\alpha}}}
\def\cb#1{\|#1\|_{C^{2,\alpha}}}

\def\lb#1{\|#1\|_{L^2}}
\def\ha#1{\|#1\|_{H^1}}
\def\hb#1{\|#1\|_{H^2}}

\def\cin{\subset\subset}
\def\Ld{\Lambda}
\def\ld{\lambda}
\def\ol{{\Omega_L}}
\def\sla{{S_L^-}}
\def\slb{{S_L^+}}
\def\e{\varepsilon}
\def\C{\mathbf{C}} %%unit cylinder
\def\cl#1{\overline{#1}}
\def\ra{\rightarrow}
\def\xra{\xrightarrow}
\def\g{\nabla}
\def\a{\alpha}
\def\b{\beta}
\def\d{\delta}
\def\th{\theta}
\def\fai{\varphi}
\def\O{\Omega}
\def\f{\frac}
\def\p{\partial}
\def\disp{\displaystyle}

\def\H{\Theta} %%truncate density function

\section{Introduction}\label{3Dintroduction}\hspace*{\parindent}

The three dimensional steady isentropic compressible Euler Equations
expressing the conservations of mass and momentum, take the following form: \be\label{da1} \left\{\ba{l}
 (\rho u_1)_{x_1}+(\rho u_2)_{x_2}+(\rho u_3)_{x_3}=0,\\
 (\rho u_1^2)_{x_1}+(\rho u_1u_2)_{x_2}+(\rho u_1u_3)_{x_3}+p_{x_1}=0,\\
 (\rho u_1u_2)_{x_1}+(\rho u_2^2)_{x_2}+(\rho u_2u_3)_{x_3}+p_{x_2}=0,\\
 (\rho u_1u_3)_{x_1}+(\rho u_2u_3)_{x_2}+(\rho u_3^2)_{x_3}+p_{x_3}=0
\ea\right. \ee where ${\bf u}=(u_1,u_2,u_3)$ is the velocity field, $p$ is
the pressure and $\rho$ is the density. We consider the ideal
polytropic gas, which means that $p(\rho)=A\rho^\gamma$, where $A$ is a positive
constant depending on the specific entropy and $\gamma \in (1,3)$ is the adiabatic
exponent. For simplicity, we take $A=\f{1}{\gamma}$. It follows from the momentum equations that the following important Bernoulli's law holds:
\begin{equation}\label{da0}
{\bf u}\cdot \nabla B=0.
\end{equation}
Here the Bernoulli's function is defined to be $B=\frac{1}{2}(u_1^2+u_2^2+u_3^2)+h(\rho)$,
where $h(\rho)$ is the enthalpy satisfying $h'(\rho)=\f{p'(\rho)}{\rho}=\f{c^2(\rho)}{\rho}$.

The original motivation to start such a research is the following: we try to construct a smooth subsonic flow in a three dimensional finitely long nozzle by formulating suitable or physically acceptable boundary conditions on upstream and downstream. Since the Euler system is
hyperbolic-elliptic coupled in a subsonic region, the main task is to understand the coupling properties of the hyperbolic and elliptic modes in subsonic flows by finding a decomposition of the Euler system to effectively decouple the different modes. We consider this as a first step to get a close look at the difficult 3-D transonic shock problems in a divergent nozzle and hope that this may shed light on that problem. For the study of transonic shock in an infinitely long nozzle, one may refer to \cite{CCF,CCS}. For the study of transonic shock problem in a duct or a De Laval nozzle, one may refer to \cite{BF,C2,C3,CY,CF,LXY1,LXY4,LXY5} and \cite{XYY,XY1,XY2,XY3} for more details. This is a natural continuation of previous results in the two dimensional case discussed in \cite{DWX,W1,W2}.

There have been many relevant studies on the 2-D steady compressible Euler equations, especially for subsonic flows around a given profile and subsonic flows in an infinitely long nozzle. Profound understanding has been achieved both physically and mathematically for subsonic flows around a given profile by Frankl and Keldysh \cite{FK}, L. Bers \cite{Bers1,Bers2}, M. Shiffman \cite{Shiffman} and many other authors \cite{Dong,FG1,FG2,Serre2,Serre3}. By a variational method, Shiffman \cite{Shiffman} proved that there is an optimal range $[0, M_{\infty})$ for Mach number $M_0$ at infinity that ensures the existence and uniqueness of a subsonic potential flow around a given profile with finite energy. Bers \cite{Bers1} presented a new proof of the existence and improved the uniqueness of a steady two-dimensional subsonic potential flow of a perfect gas around a given profile by removing the restriction of finite energy. In the recent works by Xie and Xin
\cite{XX1,XX2} for subsonic potential
flows in an infinitely long nozzle, they showed that there exists a critical value such that a global uniform subsonic flow exists uniquely in a general nozzle as long as the incoming mass flux is less than the critical value. They also obtained a class of subsonic-sonic flows by investigating the properties of these uniform subsonic flows and employing the compensated compactness method. Such a strategy
can not be applied to  the case in a multi-dimensional infinitely long nozzle ($n\geq 3$), but a
different formulation involving the potential function has been
employed by Du, Xin and Yan \cite{DXY} to obtain the existence and uniqueness
results for subsonic potential flows in a multi-dimensional infinitely long nozzle ($n\geq 2$). The new ingredients of their analysis
are methods of calculus of variations, the Moser iteration techniques for the
potential equation and a blow-up argument for infinitely long nozzles.  Concerning subsonic Euler
flows in a 2-D infinitely long nozzle, Xie and Xin \cite{XX3} explored the special structure of the 2-D Euler system and reduced the Euler system to a single second order quasi-linear elliptic equation, which  yields the existence of subsonic flows when the variation of the Bernoulli's function is sufficiently small and the mass flux is in a suitable regime with an upper critical value. For subsonic Euler flow in an axisymmetric nozzle, one may refer to \cite{DD}.

There are also some results doing this direction for both the 3-D compressible and incompressible Euler equations \cite{C3,CY,XY2,XY3,Alber,TX}. In particular, Xin and Yin \cite{XY3} developed a decomposition of the Euler equations for the uniqueness proof of the symmetric transonic shock solution and showed that the pressure satisfies a second order elliptic equation. They also studied the properties of the shock front. In \cite{C3}, the author studied the transonic shock problems in 3-D nozzle with a general section, where he gave a decomposition of Euler equations in which the elliptic part and hyperbolic part is separated at the level of principle. To deal with the loss of derivative in the velocity field, the author first obtained the vorticity by solving transport equations, then the velocity field was obtained by solving an elliptic system with complex boundary conditions. However, it seems difficult to apply these techniques to construct a subsonic Euler flow in a 3-D bounded nozzle. In \cite{Alber}, the author formulated an inflow-outflow problem for 3-D incompressible Euler system by prescribing the Bernoulli's function, the normal component of the vorticity field at the inlet of the flow region and the normal component of the velocity field on the whole boundary of the flow region. The velocity field is obtained by solving a first order elliptic system. Tang and Xin \cite{TX} have identified a class of additional boundary conditions for the vorticities and established the existence and uniqueness of solutions with non-vanishing vorticity for the three dimensional stationary incompressible Euler equations on simply connected bounded three dimensional domains with smooth boundary.

In \cite{DWX,W1,W2}, we have characterized a set of physically acceptable boundary conditions that ensure the existence and uniqueness of a subsonic irrotational flow in a finitely long flat nozzle by prescribing the incoming flow angle and the Bernoulli's function at the inlet and the end pressure at the exit. In \cite{DWX}, we show that if the incoming flow is horizontal at the inlet and an
appropriate pressure is prescribed at the exit, then there exist two positive
constants $m_0$ and $m_1$ with $m_0 < m_1$, such that a global smooth
subsonic irrotational flow exists uniquely in the nozzle, provided that the
incoming mass flux $m\in[m_0, m_1)$. The boundary conditions we have imposed in the 2-D case have a natural extension in the 3-D case. One may also prescribe the incoming flow angles ($\frac{u_2}{u_1},\frac{u_3}{u_1}$) and the Bernoulli's function $B$ at the inlet and the end pressure at the exit of the nozzle. We conjecture that such a problem should be well-posed for the 3-D Euler equations. These important clues help us a lot to find an effective decomposition for the 3-D Euler equations. It is natural to reformulate the Euler system in terms of new variables $(s=\ln \rho,\beta_2=\f{u_2}{u_1},\beta_3=\f{u_3}{u_1}, B)$. Here one should note that due to the Bernoulli's law (\ref{da0}), the Bernoulli's function $B$ will possess the same regularity as the boundary data at the inlet, which is quite different from the velocity field. So the key idea in our new formulation in terms of $(s,\beta_2,\beta_3, B)$ is to use the Bernoulli's law to reduce the dimension of the velocity field by employing $(1,\beta_2=\f{u_2}{u_1},\beta_3=\f{u_3}{u_1})$ and replace $u_1$ by $B$ through the simple algebraic formulae $u_1^2=\f{2(B-h(\rho))}{1+\beta_2^2+\beta_3^2}$. In this way, we can explore the role of the Bernoulli's law in greater depth and hope that may simplify the Euler equations a little bit. In the case of no vacuum and stagnation points, the original Euler system is equivalent to (\ref{da3}). It is shown that $s$ and $\p_2\beta_2+\p_3\beta_3$ satisfy an elliptic system. At a first sight, one may regard $(\beta_2,\beta_3)$ as hyperbolic modes since they satisfy transport equations (\ref{de7}), which will induce a loss of one derivative estimates when integrating along the particle path. Our first observation is that $W=\p_2\beta_3-\p_3\beta_2+\beta_3\p_1\beta_2-\beta_2\p_1\beta_3$ is conserved along each particle path, at least for the case of flows with a constant Bernoulli's function. Indeed, $W$ satisfies the equation (\ref{da8}), which behaves like the scaled vorticity in the 2-D case. And a physical interpretation of $W$ is also given, see section \ref{EulerConserved} for more details. This, together with the first three equations in (\ref{da3}) may help to recover the loss of one derivative in the velocity field, so that the velocity field will possess the same regularity as the pressure.

Our second observation is a system of new conservation laws, which is somehow surprising. The original motivation comes from our study on the 2-D Euler system in a curved nozzle (See \cite{W1}), where we use the flow angle $\Theta=\arctan \f{u_2}{u_1}$ instead of the angular velocity $w=\f{u_2}{u_1}$, to get a clear insight of the role played by the curvatures of the nozzle walls. We try to find a proper substitute in the 3-D case and choose the ``sine function'' of the flow angles. Finally a system of new conservation laws arises, whose physical interpretation is still unclear sofar. Even in the 2-D case, these conservation laws are new as far as we know. Yet, it seems that these information are still not enough to deal with the subsonic-sonic limit for a sequence of uniformly bounded 2-D Euler flows. How to use these new information to obtain some deeper understanding on the compressible Euler system is still under consideration.

Back to our original problem, we can construct a smooth subsonic flow in a rectangular cylinder satisfying the prescribed incoming flow angle and the Bernoulli's function at the inlet and the given end pressure at the exit, which is also required to be adjacent to some special subsonic states. We emphasize that the background can have large vorticity. We develop an iteration scheme to obtain the solution as a fixed point of a contractive operator. The velocity field $\beta_2$ and $\beta_3$ will be solved by integrating along the particle path. The corner singularities are avoided by the reflection technique. The delicate nonlinear coupling between the hyperbolic modes and elliptic modes $(\beta_2,\beta_3)$ makes it extremely difficult to have an effective iteration scheme in a 3-D nozzle with general section. How to develop an effective iteration scheme to obtain a smooth subsonic flow in a general 3-D nozzle is still under investigation.

It is worthy noting that our idea can be applied to the 3-D steady incompressible Euler equations without any essential changes. That is to say, we can employ the Bernoulli's law to get a new decomposition for the 3-D incompressible Euler system. A new conserved quantity and a system of new conservation laws can be derived as for the compressible case. As we will see, the steady incompressible Euler system is always a hyperbolic-elliptic coupled system. We can construct a smooth incompressible Euler flow adjacent to some special states in a rectangular cylinder, which will satisfy the prescribed incoming flow angles, the Bernoulli's function at the inlet and the end pressure at the exit. Our strategy can also be applied to the self-similar Euler equations, the steady Euler equations with damping, the steady Euler-Poisson equations and the steady Euler-Maxwell equations.

The rest of the paper will be organized as follows. In the next section, we will present a new formulation for the compressible Euler equations and carry out a characteristic analysis to get an insight into the structure of the Euler system. We present our main results in section \ref{EulerConserved}, where a new conserved quantity and a system of new conservation laws will be derived. Then we employ this new formulation to construct a smooth subsonic flow in a rectangular cylinder in section \ref{Subsoniceuler}. In section \ref{IEulerConserved}, we apply the same idea to give a new formulation of the incompressible Euler system. Some detailed calculations will be given in the appendix  (section \ref{appendix}).

\section{A new formulation of the 3-D Isentropic Euler equations} \label{EulerFormulation}\hspace*{\parindent}

 Since the Euler system is hyperbolic-elliptic coupled in a subsonic region, the main task is to understand the coupling properties of the hyperbolic and elliptic modes in subsonic flows by finding a decomposition of the Euler system to effectively decouple the different modes. In this section, we will give a new formulation of the Euler equations and carry out a characteristic analysis for the new system to get some insights on the coupling structure of hyperbolic and elliptic modes in the compressible Euler equations. The new reformulation has its own interest and may shed light on the difficult 3-D transonic shock problem in a de Laval nozzle. The new contribution in our formulation is to try to use the Bernoulli's law to reduce the dimension of the velocity field by defining new variables $(1,\beta_2=\frac{u_2}{u_1},\beta_3=\frac{u_3}{u_1})$ and replacing $u_1$ by the Bernoulli's function $B$ through $u_1^2=\frac{2(B-h(\rho))}{1+\beta_2^2+\beta_3^2}$.

In the whole paper, it is assumed that the flow does not contain vacuum and any
stagnation points, which means that $\rho>0$ and $(u_1,u_2,u_3)\neq(0,0,0)$ respectively.
Without loss of generality, it is assumed that $u_1>0$.

Define the following three new variables: $\beta_2=\f{u_2}{u_1},\beta_3=\f{u_3}{u_1}, s=\ln \rho$. Multiplying the first equation in (\ref{da1}) by $\f{-u_{i-1}}{\rho
u_1^2}$ , dividing the $i$-th equation in (\ref{da1}) by $\rho u_1^2$  and adding them together for $i=2,3,4$, we obtain the following new system:

\be\label{da2} \left\{\ba{l}
\p_1\beta_2+\beta_3\p_3\beta_2-\beta_2\p_3\beta_3-\beta_2\p_1s-\beta_2^2\p_2s+\f{c^2(\rho)}{u_1^2}\p_2s-\beta_2\beta_3\p_3s=0,\\
\p_1\beta_3+\beta_2\p_2\beta_3-\beta_3\p_2\beta_2-\beta_3\p_1s-\beta_2\beta_3\p_2s+\f{c^2(\rho)}{u_1^2}\p_3s-\beta_3^2\p_3s=0,\\
-\p_2\beta_2-\p_3\beta_3+\f{c^2(\rho)}{u_1^2}\p_1s-\p_1s-\beta_2\p_2s-\beta_3\p_3s=0,\\
\p_1B+\beta_2\p_2B+\beta_3\p_3B=0 \ea\right. \ee

Using the third equation, we rewrite the above system as follows:
\be\label{da3} \left\{\ba{l}
\p_1s+\beta_2\p_2s+\beta_3\p_3s-\f{c^2(\rho)}{u_1^2}\p_1s+\p_2\beta_2+\p_3\beta_3=0,\\
\p_1\beta_2+\beta_2\p_2\beta_2+\beta_3\p_3\beta_2-\f{c^2(\rho)}{u_1^2}\beta_2\p_1s+\f{c^2(\rho)}{u_1^2}\p_2s=0,\\
\p_1\beta_3+\beta_2\p_2\beta_3+\beta_3\p_3\beta_3-\f{c^2(\rho)}{u_1^2}\beta_3\p_1s+\f{c^2(\rho)}{u_1^2}\p_3s=0,\\
\p_1B+\beta_2\p_2B+\beta_3\p_3B=0. \ea\right. \ee

It is easy to verify that the above system (\ref{da3}) is equivalent
to (\ref{da1}), provided that the flow is smooth and does not contain vacuum and stagnation point.

We focus mainly on subsonic flows in which the Euler
equations become hyperbolic-elliptic coupled equations. Now
we carry out a characteristic analysis to get an insight of the structure
of the Euler equations.

We rewrite (\ref{da3}) in a matrix form as
\be\label{da4} M_1\p_1{\bf U}+M_2\p_2{\bf U}+M_3\p_3{\bf U}=0. \ee Here ${\bf
U}=(s,\beta_2,\beta_3,B)^T$ and
$$
M_1=\left[\begin{array}{rrrr}
1-\f{c^2(\rho)}{u_1^2} & 0 & 0 & 0\\
-\f{c^2(\rho)}{u_1^2}\beta_2 & 1 & 0 & 0\\
-\f{c^2(\rho)}{u_1^2}\beta_3 & 0 & 1 & 0\\
0 & 0 & 0 & 1\\
  \end{array}
\right], M_2=\left[\begin{array}{rrrr}
\beta_2 & 1 & 0 & 0\\
\f{c^2(\rho)}{u_1^2} & \beta_2 & 0 & 0\\
0 & 0 & \beta_2 & 0\\
0 & 0 & 0 & \beta_2\\
  \end{array}
\right],M_3=\left[\begin{array}{rrrr}
\beta_3 & 0 & 1 & 0\\
0 & \beta_3 & 0 & 0\\
\f{c^2(\rho)}{u_1^2} & 0 & \beta_3 & 0\\
0 & 0 & 0 & \beta_3\\
  \end{array}
\right].
$$

Denote by $\lambda$  the root of $\det(\lambda
M_1-\xi_2M_2-\xi_3M_3)$. A simple calculation shows that
$\det(\lambda
M_1-\xi_2M_2-\xi_3M_3)=(\lambda-\beta\cdot\xi)^2[(1-\f{c^2(\rho)}{u_1^2})\lambda^2-2\lambda(\beta\cdot\xi)+(\beta\cdot\xi)^2-\f{c^2(\rho)}{u_1^2}(\xi_2^2+\xi_3^2)]$,
and it has one real root $\lambda_r=\beta\cdot\xi$
with multiplicity $2$ and two conjugate complex roots
$\lambda_c^{\pm}=\lambda_R\pm i\lambda_I$, where
$\lambda_R=\f{a_2}{a_1},\lambda_I=\f{\sqrt{a_1a_3-a^2_2}}{a_1}$ and
$a_1=1-\f{c^2(\rho)}{u_1^2},a_2=\beta\cdot\xi=\beta_2\xi_2+\beta_3\xi_3,a_3=(\beta\cdot\xi)^2-\f{c^2(\rho)}{u_1^2}(\xi_2^2+\xi_3^2)$. And $a_1a_3-a^2_2=\f{c^2(\rho)}{u_1^2}[(\f{c^2(\rho)}{u_1^2}-1)(\xi_2^2+\xi_3^2)-(\beta\cdot\xi)^2]> 0$ if the flow is subsonic, i.e. $c^2(\rho)>u_1^2+u_2^2+u_3^2$.
The corresponding left eigenvectors to $\lambda_r$ are the following
$${\bf l}_r^1=(0,\xi_3+\beta_3(\beta\cdot\xi),-\xi_2-\beta_2(\beta\cdot\xi),0), {\bf l}_r^2=(0,0,0,1).$$
The corresponding left eigenvectors to $\lambda_c^{\pm}$ are
$${\bf l}_c^{\pm}={\bf l}_R\pm i{\bf l}_I=(\f{a_2}{a_1}-a_2,\xi_2,\xi_3,0)\pm i(\f{\sqrt{a_1a_3-a^2_2}}{a_1},0,0,0).$$
The differential operators corresponding to ${\bf l}_R$ and ${\bf
l}_I$ are
${\bf D_1}=(\f{c^2(\rho)}{u_1^2-c^2(\rho)}(\beta_2\p_2+\beta_3\p_3),\p_2,\p_3,0)$ and
${\bf D_2}=(1,0,0,0)$. The action of these two differential operations on (\ref{da3})
yields an elliptic system for $s$ and $\varpi=\p_2\beta_2+\p_3\beta_3$:

\be\label{da4} \left\{\ba{l}
\p_1s+\beta_2\p_2s+\beta_3\p_3s-\f{c^2(\rho)}{u_1^2}\p_1s+\varpi=0,\\
\p_1\varpi+\f{u_1^2}{u_1^2-c^2(\rho)}(\beta_2\p_2+\beta_3\p_3)\varpi+\f{c^2(\rho)}{u_1^2-c^2(\rho)}\sum_{i=2}^3\beta_i(\beta_2\p_2+\beta_3\p_3)\p_is+\\\sum_{i=2}^3\p_i(\f{c^2(\rho)}{u_1^2}\p_is)
+\f{c^2(\rho)}{u_1^2-c^2(\rho)}\p_1s(\beta_2\p_2+\beta_3\p_3)(1-\f{c^2(\rho)}{u_1^2})-\sum_{i=2}^3\p_i(\f{c^2(\rho)}{u_1^2}\beta_i)\p_1s\\
+\f{c^2(\rho)}{u_1^2-c^2(\rho)}\sum_{i=2}^3(\beta_2\p_2+\beta_3\p_3)\beta_i\p_is
+(\p_2\beta_2)^2+2\p_2\beta_3\p_3\beta_2+(\p_3\beta_3)^2=0 \ea\right.
\ee

\begin{remark}
{\it The differential operator corresponding to ${\bf l}_r^2$ is ${\bf D_4}=(0,0,0,1)$. The action of ${\bf D_4}$ on (\ref{da3}) just yields the fourth equation in (\ref{da3}), indicating that the Bernoulli's function is conserved along the particle path. One may expect that the action of the differential operator corresponding to ${\bf l}_r^1$ would show us another conserved quantity for hyperbolic modes. Indeed, the differential operator corresponding to ${\bf l}_r^1$ is
${\bf D_3}=(0,\p_3+\beta_3(\beta_2\p_2+\beta_3\p_3),-\p_2-\beta_2(\beta_2\p_2+\beta_3\p_3),0)=(0,\p_3-\beta_3\p_1+\beta_3 {\bf D},-\p_2+\beta_2\p_1-\beta_2 {\bf D},0)$, where ${\bf D}=\p_1+\beta_2\p_2+\beta_3\p_3$. Applying ${\bf D_3}$ to (\ref{da3}) and tracing the principal term, one can find that ${\bf D}W$ is the principal term involved, where $W=(\beta_3\p_1-\p_3)\beta_2-(\beta_2\p_1-\p_2)\beta_3$. However, we do not perform such a complex calculation here. It turns out that one just need to take the action of $(0,\p_3-\beta_3\p_1,-\p_2+\beta_2\p_1,0)$ on (\ref{da3}). See the next section for more details.}
\end{remark}

\begin{remark}
{\it Here $s$ and $\varpi=\p_2\beta_2+\p_3\beta_3$ satisfy an elliptic system, which is similar to the observation made in \cite{C3,CY}. Indeed, the authors in \cite{C3,CY} carried out a characteristic analysis directly to the non-isentropic Euler equations and obtained a decomposition for the Euler equations. They defined a new variable $v=\p_2 u_2+\p_3 u_3-\frac{u_2}{u_1}\p_2 u_1-\frac{u_3}{u_1}\p_3 u_1$ and showed that $v$ and the pressure $p$ satisfy the following elliptic system:
\be\label{da41} \left\{\ba{l}
\p_1 v+\lambda_1 (\p_2^2+\p_3^2)p + f_1(u_1,u_2,u_3,p,S)=0,\\
v-\lambda_2 \p_1 p=f_2(u_1,u_2,u_3,p,S).
\ea\right.\ee
Here $f_i(u_1,u_2,u_3,p,S)$ for $i=1,2$ are functions of $(u_1,u_2,u_3,p,S)$ and the first and second order derivatives of $(u_1,u_2,u_3,p,S)$.
Simple calculations show that $\varpi=\frac{v}{u_1}$. Our formulation is different from \cite{C3,CY} in the sense that we use the Bernoulli's law to reduce the dimension of the velocity field by employing $(1,\beta_2=\frac{u_2}{u_1},\beta_3=\frac{u_3}{u_1})$. Due to the Bernoulli's law (\ref{da0}), the Bernoulli's function $B$ will not lose regularity when integrating along the particle path, which is different from $u_1$. We replace $u_1$ by $B$ through the simple algebraic formulae $u_1^2=\f{2(B-h(\rho))}{1+\beta_2^2+\beta_3^2}$. }
\end{remark}

\section{A conserved quantity and a system of new conservation laws}\label{EulerConserved} \hspace*{\parindent}

 We present our main achievement in this section. First, we find a new conserved quantity $W$ satisfies an interesting transport equation, which may be used to raise the regularity of the velocity field. Second, we find a new system of conservation laws for Euler equations, some of them holds even for non-isentropic Euler equations. These new information may be used to deal with the subsonic-sonic limit for Euler flows. See Remark \ref{subsonic-sonic} for some discussions.

\subsection{A conserved quantity}\label{NewConservedquantity}\hspace*{\parindent}

Before going into the details, we state the following simple but important fact.
It is known that a two dimensional Euler flow is irrotational is equivalent to that the Bernoulli's function is a constant. However, this fact is
not true in the three dimensional case. Indeed, by employing the following identity in vector calculus
\be\label{da00}{\bf u}\cdot \nabla {\bf u}=\nabla (\frac{1}{2}|{\bf u}|^2)-{\bf u}\times curl {\bf u},\ee
and the conservation of momentum, we have
\be\label{da000}\nabla (\frac{1}{2}|{\bf u}|^2+h(\rho))-{\bf u}\times curl {\bf u}=0.\ee
So if ${\bf w}=curl {\bf u}=(w_1,w_2,w_3)\equiv 0$, then $\nabla B\equiv 0$ and $B$ is a constant. However, the constant Bernoulli's function only implies the vorticity field is parallel to the velocity field. One needs an additional condition such as $\p_2u_3-\p_3u_2\equiv 0$ to guarantee that ${\bf
w}={\bf 0}$. In a word, the constant Bernoulli's function and $\p_2u_3-\p_3u_2\equiv 0$ will imply that the flow must be
irrotational. However, $\p_2u_3-\p_3u_2$ is not invariant under Galilean transformation. Also the three component of $curl {\bf u}$ are of the same status, there is no reasons to choose one of them specially. The following interesting calculation is due to the motivation that we try to find a better substitute of $\p_2u_3-\p_3u_2$ in our new formulation. In this section, we concentrate on the flow with
a constant Bernoulli's function, in which the vorticity field is
parallel to the velocity field.

A new observation is that the quantity
$W=\p_2\beta_3-\p_3\beta_2+\beta_3\p_1\beta_2-\beta_2\p_1\beta_3=(\beta_3\p_1-\p_3)\beta_2-(\beta_2\p_1-\p_2)\beta_3$
may serve as a good candidate. That is to say, $W\equiv0$
and the constant Bernoulli's function will imply that the flow must
be irrotational. We first derive the equation satisfied by $W$ and
then give an interpretation of the physical meaning of $W$.

Since the Bernoulli's function is a constant, thus
$u_1^2=\f{2(B-h(\rho))}{1+\beta_2^2+\beta_3^2}$. To simplify
the notation, we set $G=1+\beta_2^2+\beta_3^2$ and define a new function $K(s)$ of $s$ satisfying $K'(s)=\f{c^2(e^s)}{2(B-h(e^s))}$, then we rewrite the equations
satisfied by $\beta_2$ and $\beta_3$ as follows:
\be\label{da5} \left\{\ba{l}
G^{-1}(\p_1\beta_2+\beta_2\p_2\beta_2+\beta_3\p_3\beta_2)-\beta_2\p_1K(s)+\p_2K(s)=0,\\
G^{-1}(\p_1\beta_3+\beta_2\p_2\beta_3+\beta_3\p_3\beta_3)-\beta_3\p_1K(s)+\p_3K(s)=0,\\
\ea\right. \ee

Apply $\beta_3\p_1-\p_3$ and $-\beta_2\p_1+\p_2$ to the above two
equations respectively and add them together, one can
show that $W$ satisfies the following equation.
\be\label{da6}{\bf D}W+(\p_2\beta_2+\p_3\beta_3)W-G\p_1K(s)W-G^{-1}W{\bf D}G=0.\ee

It follows from the first equation in (\ref{da3}) that $W$ satisfies
the following equation:
\be\label{da7}{\bf D}W-W {\bf D}s-G^{-1}W{\bf D} G=0.\ee
That is
\be\label{da8}{\bf D}\bigg(\f{W}{\rho G}\bigg)=0.\ee
One may refer to the appendix for the detailed computations.

Simple calculations show that
$W={\bf {\beta}}\cdot (\nabla\times {\bf \beta})=\f{u_1}{u_1}\f{w_1}{u_1}+\f{u_2}{u_1}\f{w_2}{u_1}+\f{u_3}{u_1}\f{w_3}{u_1}=\f{1}{u_1^2}({\bf u}\cdot {\bf w})$, here ${\bf \beta}=(1,\beta_2,\beta_3)$,
which is the component of ``scaled vorticity" along the ``scaled
velocity".

Then (\ref{da8}) becomes
\be\label{da81}{\bf D}\bigg(\f{{\bf u}\cdot {\bf w}}{\rho \|{\bf u}\|^2}\bigg)=0.\ee

The three components of $curl {\bf u}$ are functions dependent, so the ``scaled helicity" $\frac{W}{\rho G}$ as
a linear combination of $w_1,w_2,w_3$, which is also invariant under Galilean transformation, seems to be a better candidate for the hyperbolic mode.

Since the vorticity field is parallel to the velocity field, that is ${\bf w}=curl {\bf u}=\mu(x){\bf u}$, then we have
$$W=\f{u_1}{u_1}\f{w_1}{u_1}+\f{u_2}{u_1}\f{w_2}{u_1}+\f{u_3}{u_1}\f{w_3}{u_1}=\mu(x)(1+\beta_2^2+\beta_3^2)=\mu(x)G.$$

Then (\ref{da8}) becomes
\be\label{da82}{\bf D}\bigg(\f{\mu(x)}{\rho}\bigg)=0.\ee

Hence $\f{\mu(x)}{\rho}$ is conserved along the stream line, and plays a similar role as the ``scaled vorticity" in the two dimensional case. Recall that in 2-D case, we have $(u_1\p_1+u_2\p_2)\bigg(\f{\p_1 u_2-\p_2 u_1}{\rho}\bigg)=0$.

Historically, the stationary solution of Euler equations with the vorticity field paralleling to the velocity field was called Beltrami flow and have been investigated for over a century. Arnol'd \cite{Arnol'd} had identified a class of flows with presumably chaotic streamline, that is the Arnol'd-Beltrami-Childress (ABC) flows with ${\bf u}=(A\sin z+C \cos y, B\sin x+A\cos z,C\sin y+B \cos x)$. Turbulent flows can be understood as a superposition of Beltrami flows has been proposed in  the work of Constantin and Majda \cite{CM}. There are also some literatures on compressible Beltrami flows, for example \cite{MYZ}. Indeed, (\ref{da82}) can be derived in the following simple way: using the divergent free condition for the vorticity field, $curl {\bf u}=\mu(x){\bf u}$ will imply that $div (\mu(x){\bf u})=0$, then by employing the density equation, we directly obtain (\ref{da82}).

Now we explain the reason why we spend so much time and energy on the above complicated calculations for (\ref{da82}). The point is our calculations also work for the isentropic Euler equations, which is one of main results in this paper. Indeed, we have the following theorem.

\begin{theorem}\label{conserved}
{\it For any smooth solution $(\rho, u_1,u_2,u_3)$ to the isentropic Euler system (\ref{da1}), the quantity $W$ satisfies the following equation:
\be\begin{aligned}\label{da10}
&{\bf D}W-W{\bf D}s-G^{-1}W{\bf D}G
-\f{c^2(\rho)G}{2(B-h(\rho))^2}\bigg[(\p_2B\p_3s-\p_3B\p_2s)\\&\quad\quad\quad+\beta_2(\p_3B\p_1s-\p_1B\p_3s)+\beta_3(\p_1B\p_2s-\p_2B\p_1s)\bigg]=0.\end{aligned}\ee
That is
\be\label{da100}
{\bf D}\bigg(\f{W}{\rho G}\bigg)+\f{c^2(\rho)}{2\rho^2(B-h(\rho))^2}[(1,\beta_2,\beta_3)\cdot(\nabla \rho\times\nabla B)]=0.
\ee
}
\end{theorem}
The proof of Theorem \ref{conserved} is similar to the verification for (\ref{da8}) presented in appendix. We omit the calculations.

\begin{remark}
{\it
The equation (\ref{da100}) has its own interest, since it may be used to raise the regularity for the velocity field as explained in the following. Our original goal is to construct a smooth subsonic flow satisfying the prescribed incoming flow angles and the Bernoulli's function at the inlet and the given end pressure at the exit. Assume that the boundary data are in $C^{2,\alpha}$. Since the Bernoulli's function is conserved and will not lose regularity when integrating along the particle path, then $B \in C^{2,\alpha}$. It follows from (\ref{da3}) that $s$ satisfies a second order elliptic equation
\be\begin{aligned}\label{db20}
&\p_1\bigg((\f{c^2(\rho)}{u_1^2}-1)\p_1 s-\beta_2\p_2 s-\beta_3\p_3 s\bigg)+\p_2\bigg(-\beta_2\p_1 s+(\f{c^2(\rho)}{u_1^2}-\beta_2^2)\p_2 s-\beta_2\beta_3\p_3 s\bigg)\\&
+\p_3\bigg(-\beta_3\p_1 s-\beta_2\beta_3\p_2 s+(\f{c^2(\rho)}{u_1^2}-\beta_3^2)\p_3 s\bigg)-\bigg((\f{c^2(\rho)}{u_1^2}-1)\p_1 s-\beta_2\p_2 s-\beta_3\p_3 s\bigg)^2\\&=-\bigg((\p_2\beta_2)^2+(\p_3\beta_3)^2+2\p_2\beta_3\p_3\beta_2\bigg).
\end{aligned}\ee

This shows that $s$ also possesses $C^{2,\alpha}$ regularity in general. However, $\beta_2$ and $\beta_3$ will lose one order derivative in the process of integrating the transport equations. Roughly speaking, the quantity $W$ may help to raise the regularity of $\beta_2$ and $\beta_3$. Since $W$ satisfies (\ref{da100}), $W$ should have $C^{1,\alpha}$ regularity in general. Combining this with the first three equations in (\ref{da3}), one can obtain a divergent-curl system for $\beta_2$ and $\beta_3$, which may help to raise one more regularity for $\beta_2$ and $\beta_3$. So the velocity field would possess the same regularity as that of the pressure. The above informal argument indicates that whether the Bernoulli's function is constant or not is not so essential for enhancing the regularity of the velocity field.

}

%Roughly speaking, when one try to construct a smooth subsonic flow, the Bernoulli's function is conserved and will not lose any regularity when integrating along the particle path. It follows from (\ref{da3}) that $s$ satisfies a second order elliptic equation (see (\ref{db20}) in Remark \ref{minimumprinciple}), which possesses one more regularity than the velocity field in general. However, employing these information with (\ref{da100}), one can integrate along the particle path to obtain an integral equation for $W$, which would imply that $W$ possesses one order lower regularity than that of $s$. Combining this with the first equation in (\ref{da3}), one can obtain a first order elliptic system for $\beta_2$ and $\beta_3$, which would help to raise one more regularity for $\beta_2$ and $\beta_3$. So the velocity field would posses the same regularity as that of the pressure. The above informal argument indicates that whether the Bernoulli's function is constant or not is not an essential issue in formulating an effective iteration scheme for the construction of a smooth subsonic flow. See more discussions in the section \ref{Subsoniceuler}.

\end{remark}

\begin{remark}
{\it Indeed, one can also derive the equations (\ref{da8}) and (\ref{da81}) directly from the original Euler equations (\ref{da1}). It follows from the original Euler equations that the vorticity ${\bf w}$ satisfies the following
equations: 
\be\label{da9}
u_1\p_1{\bf w}+u_2\p_2{\bf w}+u_3\p_3{\bf
w}+div {\bf u}\cdot {\bf w}-({\bf w}\cdot\nabla) {\bf u} =0.
\ee
Combining these with the density equation, we have
\be\label{da91}(u_1\p_1+u_2\p_2+u_3\p_3)\bigg(\f{\bf w}{\rho}\bigg)-\bigg(\f{{\bf w}\cdot\nabla}{\rho}\bigg) {\bf u}=0.\ee

 It follows from (\ref{da9}) that ${\bf u}\cdot {\bf w}$ satisfies
\be\label{da92}({\bf u}\cdot\nabla)({\bf u}\cdot {\bf w}) +div {\bf u}\cdot ({\bf u}\cdot {\bf w})-({\bf w}\cdot\nabla) \|{\bf u}\|^2 =0.\ee
Employing the Bernoulli's law and the density equation, one can derive the following equation:
\be\label{da93}\|{\bf u}\|^2({\bf u}\cdot\nabla)\bigg(\f{{\bf u}\cdot {\bf w}}{\rho}\bigg) -\f{{\bf u}\cdot {\bf w}}{\rho} ({\bf u}\cdot\nabla)(\|{\bf u}\|^2)+\f{2c^2(\rho)}{\rho^2}\bigg(\|{\bf u}\|^2({\bf w}\cdot\nabla) \rho-({\bf u}\cdot {\bf w})({\bf u}\cdot\nabla) \rho\bigg) =0.\ee
That is
\be\label{da94}({\bf u}\cdot\nabla)\bigg(\f{{\bf u}\cdot {\bf w}}{\rho \|{\bf u}\|^2}\bigg)+\f{2c^2(\rho)}{\rho^2 \|{\bf u}\|^4}\bigg(\|{\bf u}\|^2({\bf w}\cdot\nabla) \rho-({\bf u}\cdot {\bf w})({\bf u}\cdot\nabla) \rho\bigg) =0.\ee

If $\nabla B\equiv {\bf 0}$, that is ${\bf u}\parallel {\bf w}$, then the last term drops out. We obtain (\ref{da81}) again. It should be emphasized that (\ref{da93}) holds also for flows with stagnation points. Hence if ${\bf u}\cdot {\bf w}\equiv 0$ at the inlet, it holds also in the whole nozzle. However, comparing (\ref{da94}) with (\ref{da100}), it seems not so clear why (\ref{da94}) can be used to raise the regularity for the velocity field.}
\end{remark}

\begin{remark}\label{R_1}{\it
We also remark here that (\ref{da000}) implies that the Bernoulli's function is conserved not only along the stream line, but also the vortex line. The last conclusion seems not getting enough attention before. Indeed, since the vorticity field is divergent free, we obtain the following new conservation law:
\be\label{da10}
div (B curl {\bf u})=0.
\ee}
\end{remark}

%\begin{remark}
%{\it Here it should be emphasized that the key idea in the above calculations is to use the Bernoulli's law to reduce the dimension of the velocity field by employing $(1,\beta_2=\f{u_2}{u_1},\beta_3=\f{u_3}{u_1})$ and $u_1^2=\f{2(B-h(\rho))}{1+\beta_2^2+\beta_3^2}$. Thus the Bernoulli's law plays a basic role in our new formulation and detailed calculations.}
%\end{remark}

\subsection {A system of new conservation laws} \label{NewConservationlaws}\hspace*{\parindent}

In the following, we assume that the Bernoulli's function is a
constant and the speed $|{\bf u}|$ has a positive lower bound $\delta$. It is easy to see that $G$ satisfies the following
Riccati-type equation:
\be\label{da11}{\bf D}G-\f{c^2(\rho)}{(B-h(\rho))}\p_1s
G^2+\f{c^2(\rho)}{(B-h(\rho))}G{\bf D}s=0.\ee
Define three new variables:
$K_1=G^{-\f{1}{2}},K_i=\f{\beta_i}{G^{\f{1}{2}}}=\beta_iK_1,i=2,3$.
Here $K_2$ and $K_3$ behave like the ``sine function" of ``the flow
angles". In the following, we derive the equations satisfied
by $K_2$ and $K_3$. To simplify the notation, we define two new
functions $I(\rho)=e^{\int_0^{\rho}\f{c^2(t)}{2t(B-h(t))}dt}$ and
$Q(\rho)=\int_0^{\rho}\f{c^2(t)}{2(B-h(t))}I(t)^{-2}dt$. Since the speed $|{\bf u}|$ has a lower bound  $\delta$, $\int_0^{\rho}\f{c^2(t)}{2t(B-h(t))}dt\leq \frac{1}{\delta^2} h(\rho)$, so $I(\rho)$ and $Q(\rho)$ are well-defined and $Q(\rho)>0$ for $\rho>0$. One should note that
$I(\rho)$ satisfies
$\f{I'(\rho)}{I(\rho)}=\f{c^2(\rho)}{2\rho(B-h(\rho))}$.

It follows from (\ref{da11}) that $K_1$ satisfies:
\be\label{da12}
{\bf D}K_1-\f{c^2(\rho)}{2\rho(B-h(\rho))}K_1 {\bf D}\rho+\f{c^2(\rho)}{2\rho(B-h(\rho))}K_1^{-1}\p_1\rho=0.\ee
While the first equation in (\ref{da3}) yields
\be\label{da120}
\f{K_1}{\rho}{\bf D}\rho-\f{c^2(\rho)}{2\rho(B-h(\rho))}K_1^{-1} \p_1\rho+K_1(\p_2\beta_2+\p_3\beta_3)=0.\ee
Combing these two equations, we obtain
\be\label{da121}
\f{K_1}{\rho}{\bf D}\rho+{\bf D}K_1-\f{c^2(\rho)}{2\rho(B-h(\rho))}K_1 {\bf D}\rho+K_1(\p_2\beta_2+\p_3\beta_3)=0.\ee
That is
\be\label{da122}
\f{1}{\rho}{\bf K}\cdot\nabla\rho+div {\bf K}-\f{c^2(\rho)}{2\rho(B-h(\rho))}{\bf K}\cdot\nabla\rho=0.\ee
Using the definition of $I(\rho)$, we have the following new conservation law, which is similar to the conservation of mass:
\be\label{da123}
div (\rho {\bf A})=0.
\ee
 where ${\bf A}=\f{1}{I(\rho)}{\bf K}=\f{1}{I(\rho)}(K_1,K_2,K_3)$.

Similarly, one can also derive the equations satisfied by $K_i, i=1,2,3$, as
\be\label{da13} \left\{\ba{l}
{\bf D}K_1-\f{c^2(\rho)}{2\rho(B-h(\rho))}K_1{\bf D}\rho+\f{c^2(\rho)}{2\rho(B-h(\rho))}K_1^{-1}\p_1\rho=0.\\
{\bf D}K_2-\f{c^2(\rho)}{2\rho(B-h(\rho))}K_2{\bf D}\rho+\f{c^2(\rho)}{2\rho(B-h(\rho))}K_1^{-1}\p_2\rho=0,\\
{\bf D}K_3-\f{c^2(\rho)}{2\rho(B-h(\rho))}K_3{\bf D}\rho+\f{c^2(\rho)}{2\rho(B-h(\rho))}K_1^{-1}\p_3\rho=0.
\ea\right. \ee
That is:
\be\label{da131} \left\{\ba{l}
{\rho\bf K}\cdot \nabla K_1-K_1\f{\rho{\bf K}\cdot\nabla I(\rho)}{I(\rho)} +\f{c^2(\rho)}{2(B-h(\rho))}\p_1\rho=0.\\
{\rho\bf K}\cdot \nabla K_2-K_2\f{\rho{\bf K}\cdot\nabla I(\rho)}{I(\rho)}+\f{c^2(\rho)}{2(B-h(\rho))}\p_2\rho=0,\\
{\rho\bf K}\cdot \nabla K_3-K_3\f{\rho{\bf K}\cdot\nabla I(\rho)}{I(\rho)}+\f{c^2(\rho)}{2(B-h(\rho))}\p_3\rho=0.
\ea\right. \ee

With the help of (\ref{da123}), we
obtain the following conservation laws, which are similar to the conservation of momentum:
\be\label{da15} \left\{\ba{l}
\p_1(\rho A_1^2)+\p_2(\rho A_1A_2)+\p_3(\rho A_1A_3)+\p_1 Q(\rho)=0,\\
\p_1(\rho A_1A_2)+\p_2(\rho A_2^2)+\p_3(\rho A_2A_3)+\p_2 Q(\rho)=0,\\
\p_1(\rho A_1A_3)+\p_2(\rho A_2A_3)+\p_3(\rho A_2^3)+\p_3 Q(\rho)=0.
\ea\right. \ee

For the isentropic Euler equations, since $(\p_1+\beta_2\p_2+\beta_3\p_3)B=0$, we still have
\be\label{da101}
\f{1}{\rho}{\bf K}\cdot\nabla\rho+div {\bf K}-\f{1}{I(\rho,B)}{\bf K}\cdot\nabla I(\rho,B)=0,\ee
where $I(\rho,B)=e^{\int_0^{\rho}\f{c^2(t)}{2t(B-h(t))}dt}$. Hence one can still obtain the following conservation law:
\be\label{da102}
\p_1\bigg(\rho \f{K_1}{I(\rho, B)}\bigg)+\p_2 \bigg(\rho \f{K_2}{I(\rho, B)}\bigg)+\p_3\bigg(\rho \f{K_3}{I(\rho, B)}\bigg)=0.\ee
However, similar conservation laws as (\ref{da15}) are not true in general.

For the non-isentropic Euler equations, one can perform the same operations to get similar information. Actually, the non-isentropic Euler equations can be rewritten as
\be\label{db10}\left\{\ba{l}
\f{1}{\rho}(\p_1\rho+\beta_2\p_2\rho+\beta_3\p_3\rho)-\f{1}{\rho u_1^2}\p_1P+\p_2\beta_2+\p_3\beta_3=0,\\
\p_1\beta_2+\beta_2\p_2\beta_2+\beta_3\p_3\beta_2-\f{1}{\rho u_1^2}\beta_2\p_1P+\f{1}{\rho u_1^2}\p_2P=0,\\
\p_1\beta_3+\beta_2\p_2\beta_3+\beta_3\p_3\beta_3-\f{1}{\rho u_1^2}\beta_3\p_1P+\f{1}{\rho u_1^2}\p_3P=0,\\
\p_1B+\beta_2\p_2B+\beta_3\p_3B=0,\\
\p_1S+\beta_2\p_2S+\beta_3\p_3S=0.
\ea\right. \ee

In this case, $B=\f{1}{2}(u_1^2+u_2^2+u_3^2)+h(\rho,S)=\f{1}{2}u_1^2(1+\beta_2^2+\beta_3^2)+h(\rho,S)$, where $h(\rho,S)=\int_0^{\rho}\frac{\p_{\rho}P(\tau,S)}{\tau}d{\tau}$. Hence $u_1^2=\f{2(B-h(\rho,S))}{1+\beta_2^2+\beta_3^2}$.

Similar calculations as previous sections show that $W=\p_2\beta_3-\p_3\beta_2+\beta_3\p_1\beta_2-\beta_2\p_1\beta_3$ satisfies the following equation:
\be\begin{aligned}\label{db12}
&{\bf D}\bigg(\f{W}{\rho G}\bigg)+\f{1}{2\rho^3(P,S)(B-h(P,S))^2}\bigg\{(B-h(P,S))\p_S\rho(P,S){\bf {\beta}}\cdot(\nabla P\times \nabla S)\\&
\quad\quad\quad \quad +\rho(P,S)[{\bf \beta}\cdot(\nabla P\times \nabla B)-\p_Sh(P,S) {\bf \beta}\cdot(\nabla P\times \nabla S)]\bigg\}=0,
\end{aligned}\ee
where we have rewritten $\rho$ and $h$ as smooth functions of $(P,S)$.

Also a similar conservation law holds:
\be\label{db16}
\p_1\bigg(\rho \f{K_1}{I(\rho, B, S)}\bigg)+\p_2 \bigg(\rho \f{K_2}{I(\rho, B, S)}\bigg)+\p_3\bigg(\rho \f{K_3}{I(\rho, B, S)}\bigg)=0,
\ee
where $I(\rho,B,S)=e^{\int_0^{\rho}\f{\p_{\rho} P(t,S)}{2t(B-h(t, S))}dt}$.

We summarize the above results in the following theorem, which is our second main results in this paper.
\begin{theorem}\label{law}
{\it For any smooth solution $(\rho, u_1,u_2,u_3)$ to the isentropic Euler equations (\ref{da1}), if the Bernoulli's function $B$ is a constant throughout the flow regime, we have the following new conservation laws
\be\label{law1}
\begin{aligned}
& div (\rho {\bf A})=0,\\
& \p_1(\rho A_1^2)+\p_2(\rho A_1A_2)+\p_3(\rho A_1A_3)+\p_1 Q(\rho)=0,\\
& \p_1(\rho A_1A_2)+\p_2(\rho A_2^2)+\p_3(\rho A_2A_3)+\p_2 Q(\rho)=0,\\
& \p_1(\rho A_1A_3)+\p_2(\rho A_2A_3)+\p_3(\rho A_2^3)+\p_3 Q(\rho)=0.
\end{aligned}
\ee
For the general isentropic Euler equations, we also have
\be\label{law2}
\p_1\bigg(\rho \f{K_1}{I(\rho, B)}\bigg)+\p_2 \bigg(\rho \f{K_2}{I(\rho, B)}\bigg)+\p_3\bigg(\rho \f{K_3}{I(\rho, B)}\bigg)=0.
\ee
Finally, for the non-isentropic Euler equations, we also have
\be\label{law3}
\p_1\bigg(\rho \f{K_1}{I(\rho, B, S)}\bigg)+\p_2 \bigg(\rho \f{K_2}{I(\rho, B, S)}\bigg)+\p_3\bigg(\rho \f{K_3}{I(\rho, B, S)}\bigg)=0.
\ee
}
\end{theorem}

\begin{remark}{\it
Since $G$ satisfies a Riccati-type equation (\ref{da11}), one may expect some blow-up results. However, the ambiguous sign of $\p_1s$ in the coefficient of $G^2$ in (\ref{da11}) makes the whole argument nontrivial. One should note that $G$ blows up means that the fluid will turn around in the flow region.}
\end{remark}

\begin{remark}
{\it
The physical interpretations of (\ref{da123}) and (\ref{da15}) are not clear at this moment. How to use these new information to obtain some useful results is an interesting problem and will be investigated further.}
\end{remark}

\begin{remark}{\it
One may perform the same operation as we have done in Remark \ref{R_1}, to (\ref{da123}) and (\ref{da15}) to obtain a new conservation law involving $curl {\bf A}$. However, one should note that the Bernoulli's function for (\ref{da15}) is just a constant, so no new information can be obtained in this case.}
\end{remark}

\begin{remark}{\it
One may perform similar operations as we have done for the original Euler equations, to the new conservation laws (\ref{da123}) and (\ref{da15}), it turns out that one will get (\ref{da123}) and (\ref{da15}) again. That is, (\ref{da123}) and (\ref{da15}) is invariant under our operations.}
\end{remark}

%\begin{remark}{\it} \end{remark}

%\begin{remark}{\it}\end{remark}

\begin{remark}{\it
These calculations certainly work for the two dimensional Euler system. Indeed, $W$ will always be zero. We can obtain a system of new conservation laws as (\ref{da123}) and (\ref{da15}) for potential flows and a new conservation law as (\ref{da102}) or (\ref{db16}) for the 2-D Euler equations. Even in this case, all these information are completely new as far as we know. %However, it seems that these new information are still not enough to deal with the subsonic-sonic limit for the 2-D Euler flow.

}
\end{remark}

\begin{remark}\label{subsonic-sonic}
{\it
As shown in \cite{cdsw} and \cite{XX1,XX2}, a 2-D potential sonic flow was obtained by compensated compactness from a sequence of smooth subsonic potential flows. A natural interesting question is whether one can show the existence of Euler sonic flows through similar subsonic-sonic limit argument. The new conservation laws obtained here may help, however, a rigourous mathematical proof is still not available.  
}
\end{remark}

\begin{remark}{\it
We may apply the same idea to the self-similar compressible Euler equations, which take the following form:
\be\label{db17} \left\{\ba{l}
 (\rho U_1)_{\xi_1}+(\rho U_2)_{\xi_2}+(\rho U_3)_{\xi_3}+3\rho=0,\\
 \rho U_1 \p_{\xi_1}U_1+\rho U_2 \p_{\xi_2}U_1+\rho U_3 \p_{\xi_3}U_1+\p_{\xi_1}p+\rho U_1=0,\\
 \rho U_1 \p_{\xi_1}U_2+\rho U_2 \p_{\xi_2}U_2+\rho U_3 \p_{\xi_3}U_2+\p_{\xi_2}p+\rho U_2=0,\\
 \rho U_1 \p_{\xi_1}U_3+\rho U_2 \p_{\xi_2}U_3+\rho U_3 \p_{\xi_3}U_3+\p_{\xi_3}p+\rho U_3=0.
\ea\right.\ee

Here ${\bf \xi}=(\xi_1,\xi_2, \xi_3)=(\f{x_1}{t},\f{x_2}{t},\f{x_3}{t})$ and $(U_1,U_2,U_3)({\bf \xi})=(u_1,u_2,u_3)({\bf \xi})-(\xi_1,\xi_2, \xi_3)$. The main difference in (\ref{db17}) is that the Bernoulli's law does not hold any more. One can still define $s=\ln \rho, \beta_2=\f{U_2}{U_1}, \beta_3=\f{U_3}{U_1},B=\f{1}{2}(U_1^2+U_2^2+U_3^2)+h(\rho)$ and obtain a new formulation for (\ref{db17}):
\be\label{db18}\left\{\ba{l}
(\p_{\xi_1}s+\beta_2\p_{\xi_2}s+\beta_3\p_{\xi_3}s)-\f{c^2(\rho)}{ U_1^2}\p_{\xi_1}s+\p_2\beta_2+\p_3\beta_3+\f{2}{U_1}=0,\\
\p_{\xi_1}\beta_2+\beta_2\p_{\xi_2}\beta_2+\beta_3\p_{\xi_3}\beta_2-\f{c^2(\rho)}{ U_1^2}\beta_2\p_{\xi_1}s+\f{c^2(\rho)}{U_1^2}\p_{\xi_2}s=0,\\
\p_{\xi_1}\beta_3+\beta_2\p_{\xi_2}\beta_3+\beta_3\p_{\xi_3}\beta_3-\f{c^2(\rho)}{U_1^2}\beta_3\p_{\xi_1}s+\f{c^2(\rho)}{U_1^2}\p_{\xi_3}s=0,\\
\p_{\xi_1}B+\beta_2\p_{\xi_2}B+\beta_3\p_{\xi_3}B+\f{(U_1^2+U_2^2+U_3^2)}{U_1}=0.
\ea\right. \ee
Although the Bernoulli's law does not hold in this case, the Bernoulli's function $B$ still has the same regularity as that of the pressure, and $W=\p_2\beta_3-\p_3\beta_2+\beta_3\p_1\beta_2-\beta_2\p_1\beta_3$  may also help to raise the regularity of $\beta_2$ and $\beta_3$. Hence the reformulation (\ref{db18}) may provide some new insights on the self-similar Euler equations. Same arguments work for the steady Euler equations with damping.}
\end{remark}

\begin{remark}{\it
In \cite{W3}, we will employ the same idea to investigate the structural stability of some steady subsonic solutions for 1-D Euler-Poisson system under multi-dimensional perturbations of suitable boundary conditions and background charge.}
\end{remark}

\section{Subsonic Euler flows in a rectangular cylinder}\label{Subsoniceuler} \hspace*{\parindent}

In this section, we try to construct a subsonic Euler flow in a rectangular cylinder by imposing suitable boundary conditions at the inlet and exit. The cylinder will be $\O=[0,1]\times[0,1]\times[0,1]$. Actually the flow we construct below will be adjacent to some special subsonic states $(\rho_0,u_0(x_2,x_3),0,0)$. Here $\rho_0$ is a constant, $u_0(x_2,x_3)\in C^{3,\alpha}([0,1]\times[0,1])$ satisfying $0<u_0(x_2,x_3)<c(\rho_0)$. We also set $B_0(x_2,x_3)=\f{1}{2}u_0^2(x_2,x_3)+h(\rho_0), s_0=\ln\rho_0$.

At the inlet of the nozzle $x_1=0$, we impose the flow angles and the Bernoulli's function:
\be\label{de1}\left\{\ba{l}
\beta_i(0,x_2,x_3)=\epsilon \beta_{i}^{in}(x_2,x_3), i=2,3,\\
B(0,x_2,x_3)= B_0(x_2,x_3)+\epsilon B^{in}(x_2,x_3).
\ea\right.\ee
Here the compatibility conditions should be satisfied:
\be\label{de100}\left\{\ba{l}
\p_2^j\beta_{2}^{in}(0,x_3)=\p_2^j\beta_{2}^{in}(1,x_3)=0, \\
\p_3^j\beta_{3}^{in}(x_2,0)=\p_3^j\beta_{3}^{in}(x_2,1)=0,\ \ j=0,2,\\
\p_2^k B^{in}(0,x_3)=\p_2^k B^{in}(1,x_3)=\p_3^k B^{in}(x_2,0)=\p_3^k B^{in}(x_2,1)=0,\ \ k=1,3.
\ea\right.\ee

At the exit of the nozzle $x_1=1$, we prescribe the end pressure:
\be\label{de2}
p(1,x_2,x_3)=\f{1}{\gamma} e^{\gamma(s_0+\epsilon s_e(x_2,x_3))}.
\ee
Here we also require that $s_e$ satisfies the following compatibility conditions:
\be\label{de101}\left\{\ba{l}
\p_2^j s_e(0,x_3)=\p_2^j s_e(1,x_3)=0,\\
\p_3^j s_e(x_2,0)=\p_3^j s_e(x_2,1)=0, \ \ j=1,3.
\ea\right.\ee

While on the nozzle walls, the usual slip boundary condition is imposed:
\be\label{de3}\left\{\ba{l}
u_2(x_1,0,x_3)=u_2(x_1,1,x_3)=0,\\
u_3(x_1,x_2,0)=u_3(x_1,x_2,1)=0.
\ea\right.\ee

Mathematically, we are going to prove that (\ref{da1}) with boundary conditions (\ref{de1}),(\ref{de2}) and (\ref{de3}) satisfying compatibility conditions (\ref{de100}) and(\ref{de101}), has a unique subsonic solution.
%\be\label{e4} \left\{\ba{l}
% (\rho u_1)_{x_1}+(\rho u_2)_{x_2}+(\rho u_3)_{x_3}=0,\\
% (\rho u_1^2)_{x_1}+(\rho u_1u_2)_{x_2}+(\rho u_1u_3)_{x_3}+p_{x_1}=0,\\
% (\rho u_1u_2)_{x_1}+(\rho u_2^2)_{x_2}+(\rho u_2u_3)_{x_3}+p_{x_2}=0,\\
% (\rho u_1u_3)_{x_1}+(\rho u_2u_3)_{x_2}+(\rho u_3^2)_{x_3}+p_{x_3}=0.\\
% \f{u_i}{u_1}(0,x_2,x_3)=\epsilon \beta_{i0}(x_2,x_3),\ \ \text{for}\ \ i=2,3,  \\
% B(0,x_2,x_3)=B_c+\epsilon B_0(x_2,x_3), \  \text{and} \ \ u_1(0,x_2,x_3) > 0,\\
% u_2(x_1,0,x_3)=u_2(x_1,1,x_3)=0,\\
%u_3(x_1,x_2,0)=u_3(x_1,x_2,1)=0,\\
%p(1,x_2,x_3)=p_c+\epsilon p_e(x_2,x_3).
%\ea\right. \ee

\subsection{Extension to the domain $\O_e=[0,1]\times \mathrm{T}^2$}\label{Eulerextension}\hspace*{\parindent}

Suppose that the flow $(\rho, u_1,u_2,u_3)\in C^{3,\alpha}(\bar \O)\times C^{2,\alpha}(\bar \O)^3$ to be constructed has the following properties:
\be\label{de41} \left\{\ba{l}
\p_2^j(\rho,u_1)(x_1,0,x_3)=\p_2^j(\rho,u_1)(x_1,1,x_3)=0,\\
\p_3^j(\rho,u_1)(x_1,x_2,0)=\p_3^j(\rho,u_1)(x_1,x_2,1)=0,\ \ j=1,3,\\
\p_2^ku_2(x_1,0,x_3)=\p_2^ku_2(x_1,1,x_3)=0,\\
\p_3^ku_3(x_1,x_2,0)=\p_3^ku_3(x_1,x_2,1)=0,\ \ k=0,2.
\ea\right. \ee
Then we may extend $(\rho, u_1,u_2,u_3)$ in the following way, to $(\hat\rho, \hat u_1,\hat u_2,\hat u_3)\in C^{3,\alpha}([0,1]\times \mathbb{R}^2)\times C^{2,\alpha}([0,1]\times \mathbb{R}^2)^3$:

For $(x_1,x_2,x_3)\in [0,1]\times [-1,1]\times[-1,1]$, we define $(\hat\rho,\hat u_1,\hat u_2,\hat u_3)$ as follows
$$
(\hat\rho,\hat u_1,\hat u_2,\hat u_3)({\bf x})=\left\{
    \ba{ll}
    (\rho,u_1, u_2, u_3)(x_1,x_2,x_3),\ \ \ \ \ \text{ if } (x_2,x_3)\in [0,1]\times [0,1],\\
    (\rho,u_1, -u_2, u_3)(x_1,-x_2,x_3),\ \ \ \ \text{ if } (x_2,x_3)\in [-1,0]\times [0,1],\\
    (\rho,u_1, u_2, -u_3)(x_1,x_2,-x_3),\ \ \ \ \text{ if } (x_2,x_3)\in [0,1]\times [-1,0],\\
    (\rho,u_1, -u_2, -u_3)(x_1,-x_2,-x_3),\ \ \ \text{ if } (x_2,x_3)\in [-1,0]\times [-1,0].
    \ea
    \right.
$$
Then we extend $(\hat\rho, \hat u_1,\hat u_2,\hat u_3)$ periodically to $[0,1]\times \mathbb{R}^2$ with period $2$. It is easy to verify that $(\hat\rho, \hat u_1,\hat u_2,\hat u_3)$ belong to $C^{3,\alpha}([0,1]\times \mathbb{R}^2)\times C^{2,\alpha}([0,1]\times \mathbb{R}^2)^3$. Moreover, $(\hat\rho, \hat u_1,\hat u_2,\hat u_3)$ satisfy the Compressible Euler equations. Due to these reasons, one may directly work on the domain $[0,1]\times \mathrm{T}^2$ (Here $\mathrm{T}^2$ is a 2-torus), so that the difficulty caused by corner singularities and slip boundary conditions can be avoided. We may extend $(u_0,B_0,\beta_{2}^{in},\beta_{3}^{in},B^{in}, s_e)$ to $\mathrm{T}^2$, which will still be denoted by $(u_0,B_0,\beta_{2}^{in},\beta_{3}^{in},B^{in}, s_e)$.

\subsection{Main results}\label{3dEulermainresults}\hspace*{\parindent}

The main result in this section is the following existence and uniqueness theorem.
\begin{theorem}\label{ETH}
Given $(\beta_{2}^{in},\beta_{3}^{in},B^{in}, s_e)\in C^{3,\alpha}(\mathrm{T}^2)$ satisfying the compatibility conditions (\ref{de100}) and (\ref{de101}), there exists a positive small number $\epsilon_0$, which depends on the background subsonic state $(\rho_0,u_0(x_2,x_3),0,0)$ and $(\beta_{2}^{in},\beta_{3}^{in},B^{in}, s_e)$, such that if $0<\epsilon<\epsilon_0$, then there exists a unique smooth subsonic flow $(\rho,u_1,u_2,u_3)$$\in C^{3,\alpha}(\O_e)\times (C^{2,\alpha}(\O_e))^3$ to (\ref{da1}) satisfying boundary conditions (\ref{de1}),(\ref{de2}) and (\ref{de3}). Moreover, the following estimate holds:
\be\label{de50}
\|(u_1,u_2,u_3)-(u_0(x_2,x_3),0,0)\|_{C^{2,\alpha}(\O_e)}+\|\rho-\rho_0\|_{C^{3,\alpha}(\O_e)}\leq C\epsilon.
\ee
Here $C$ is a constant depending on $(\rho_0,u_0(x_2,x_3),0,0)$ and $(\beta_{2}^{in},\beta_{3}^{in},B^{in}, s_e)$.
\end{theorem}

Since we are looking for a smooth subsonic Euler flow in the whole nozzle, which is also close to some special subsonic states, no vacuum and stagnation points will appear in the flow region, we may use the new formulation developed in the section \ref{EulerFormulation} and section \ref{EulerConserved}. Hence to prove the above theorem, it is equivalent to show that the following system has a unique smooth subsonic solution:
\be\label{de5} \left\{\begin{array}{ll}
\p_1s+\beta_2\p_2s+\beta_3\p_3s-\f{c^2(\rho)}{u_1^2}\p_1s+\p_2\beta_2+\p_3\beta_3=0,\\
\p_1\beta_2+\beta_2\p_2\beta_2+\beta_3\p_3\beta_2-\f{c^2(\rho)}{u_1^2}\beta_2\p_1s+\f{c^2(\rho)}{u_1^2}\p_2s=0,\\
\p_1\beta_3+\beta_2\p_2\beta_3+\beta_3\p_3\beta_3-\f{c^2(\rho)}{u_1^2}\beta_3\p_1s+\f{c^2(\rho)}{u_1^2}\p_3s=0,\\
\p_1B+\beta_2\p_2B+\beta_3\p_3B=0,\\
\beta_i(0,x_2,x_3)=\epsilon \beta_{i}^{in}(x_2,x_3),\ \ \text{for}\ \ i=2,3,  \\
B(0,x_2,x_3)=B_0+\epsilon B^{in}(x_2,x_3), \  \text{and} \ \ u_1(0,x_2,x_3)> 0,  \\
s(1,x_2,x_3)=s_0+\epsilon s_e(x_2,x_3).
\end{array}\right. \ee
We can derive that $\bar s=s-s_0$, $\beta_2,\beta_3$ and $\bar B=B-B_0$ satisfy the following equations, respectively:
\be\label{de6} \left\{\begin{array}{ll}
\p_1\bigg((\f{c^2(\rho)}{u_1^2}-1)\p_1 \bar s-\beta_2\p_2\bar s-\beta_3\p_3\bar s\bigg)+\p_2\bigg(-\beta_2\p_1\bar s+(\f{c^2(\rho)}{u_1^2}-\beta_2^2)\p_2\bar s-\beta_2\beta_3\p_3\bar s\bigg)\\
+\p_3\bigg(-\beta_3\p_1\bar s-\beta_2\beta_3\p_2\bar s+(\f{c^2(\rho)}{u_1^2}-\beta_3^2)\p_3\bar s\bigg)-\bigg((\f{c^2(\rho)}{u_1^2}-1)\p_1 \bar s-\beta_2\p_2\bar s-\beta_3\p_3\bar s\bigg)^2\\+(\p_2\beta_2)^2+(\p_3\beta_3)^2+2\p_2\beta_3\p_3\beta_2=0,\\
-\bigg((\f{c^2}{u_1^2}-1)\p_1 \bar s-\beta_2\p_2\bar s-\beta_3\p_3\bar s\bigg)(0,x_2,x_3)=\epsilon(\p_2\beta_{2}^{in}+\p_3\beta_{3}^{in}),\\
\bar s(1,x_2,x_3)=\epsilon s_e(x_2,x_3).
\end{array}\right. \ee

\be\label{de7} \left\{\begin{array}{ll}
\p_1\beta_2+\beta_2\p_2\beta_2+\beta_3\p_3\beta_2-\f{c^2(\rho)}{u_1^2}\beta_2\p_1\bar s+\f{c^2(\rho)}{u_1^2}\p_2\bar s=0,\\
\p_1\beta_3+\beta_2\p_2\beta_3+\beta_3\p_3\beta_3-\f{c^2(\rho)}{u_1^2}\beta_3\p_1\bar s+\f{c^2(\rho)}{u_1^2}\p_3\bar s=0,\\
\beta_2(0,x_2,x_3)=\epsilon \beta_{2}^{in}(x_2,x_3),\\
\beta_3(0,x_2,x_3)=\epsilon \beta_{3}^{in}(x_2,x_3),\\
\end{array}\right. \ee

\be\label{de8} \left\{\begin{array}{ll}
\p_1\bar B+\beta_2\p_2\bar B+\beta_3\p_3\bar B=-\beta_2\p_2B_0-\beta_3\p_3B_0,\\
 \bar B(0,x_2,x_3)=\epsilon B^{in}(x_2,x_3).
\end{array}\right. \ee

\subsection{Proof of Theorem \ref{ETH}}\label{3dEulerproof}\hspace*{\parindent}

The main idea is simple: we construct an operator $\Lambda$ on a suitable space, which will be bounded in a high order norm and contraction in a low order norm.

The solution class will be given by
$$\Xi=\{(\bar s,\beta_2,\beta_3,\bar B): \|\bar s\|_{C^{3,\alpha}(\O_e)}\leq \delta, \|\beta_2,\beta_3,\bar B\|_{C^{2,\alpha}(\O_e)}\leq \delta.\}$$
Here $\delta$ will be determined later. For any given $(\t{\bar s},\t \beta_2,\t \beta_3,\t {\bar B})\in \Xi$, we define an operator $\Lambda:(\t{\bar s},\t \beta_2,\t \beta_3,\t {\bar B})\mapsto (\bar s, \beta_2, \beta_3,\bar B)$ mapping from $\Xi$ to itself, through the following iteration scheme.

Step 1. To obtain $\bar s$ by solving the following linearized elliptic system:
\be\label{de9} \left\{\begin{array}{ll}
\p_1\bigg((\f{c^2(\t \rho)}{{\t u_1}^2}-1)\p_1 \bar s-\t\beta_2\p_2\bar s-\t\beta_3\p_3\bar s\bigg)+\p_2\bigg(-\t\beta_2\p_1\bar s+(\f{c^2(\t \rho)}{{\t u_1}^2}-{\t\beta_2}^2)\p_2\bar s-\t\beta_2\t\beta_3\p_3\bar s\bigg)\\ \quad \quad\quad
+\p_3\bigg(-\t\beta_3\p_1\bar s-\t\beta_2\t\beta_3\p_2\bar s+(\f{c^2(\t \rho)}{{\t u_1}^2}-{\t\beta_3}^2)\p_3\bar s\bigg)=F_1,\\
-\bigg((\f{c^2(\t \rho)}{{\t u_1}^2}-1)\p_1 \bar s-\t\beta_2\p_2\bar s-\t\beta_3\p_3\bar s\bigg)(0,x_2,x_3)=\epsilon(\p_2\beta_{2}^{in}+\p_3\beta_{3}^{in}),\\
\bar s(1,x_2,x_3)=\epsilon s_e(x_2,x_3).
\end{array}\right. \ee
Here $F_1=\bigg((\f{ c^2(\t \rho)}{{\t u_1}^2}-1)\p_1 \t{\bar s}-\beta_2\p_2\t{\bar s}-\beta_3\p_3\t{\bar s}\bigg)^2-\bigg((\p_2\t\beta_2)^2+(\p_3\t\beta_3)^2+2\p_2\t\beta_3\p_3\t\beta_2\bigg)$, $\t \rho=e^{\t {\bar s}+s_0}$ and ${\t u_1}^2=\frac{2(B_0+\t{\bar B}-h(\t \rho))}{1+\t\beta_2^2+\t\beta_3^2}$.

By the standard elliptic estimates, we can show that (\ref{de9}) has a unique solution $\bar s\in C^{3,\alpha}(\bar\O)$ with the following estimate:
\be\begin{aligned}\label{de10}
\|\bar s\|_{C^{3,\alpha}(\O_e)}&\leq C_1\bigg(\|F_1\|_{C^{1,\alpha}(\O_e)}+\epsilon\|\p_2\beta_{2}^{in}+\p_3\beta_{3}^{in}\|_{C^{2,\alpha}(\mathrm{T}^2)}+\epsilon\|s_e\|_{C^{3,\alpha}(\mathrm{T}^2)}\bigg)\\&\leq C_2(\delta^2+\epsilon).
\end{aligned}\ee
Here $C_1$ depends only on the background solution and $\O_e$, and $C_2$ depends also on $(\beta_{2}^{in},\beta_{3}^{in},B^{in}, s_e)$.

Step 2. To obtain $\beta_2,\beta_3$ by solving the following hyperbolic equations:
\be\label{de11} \left\{\begin{array}{ll}
\p_1\beta_2+\t\beta_2\p_2\beta_2+\t\beta_3\p_3\beta_2-\f{c^2(\t\rho)}{{\t u_1}^2}\t\beta_2\p_1\bar s+\f{c^2(\t\rho)}{{\t u_1}^2}\p_2\bar s=0,\\
\p_1\beta_3+\t\beta_2\p_2\beta_3+\t\beta_3\p_3\beta_3-\f{c^2(\t\rho)}{{\t u_1}^2}\t\beta_3\p_1\bar s+\f{c^2(\t\rho)}{{\t u_1}^2}\p_3\bar s=0,\\
\beta_2(0,x_2,x_3)=\epsilon\beta_{2}^{in}(x_2,x_3),\\
\beta_3(0,x_2,x_3)=\epsilon\beta_{3}^{in}(x_2,x_3).
\end{array}\right. \ee

The particle path $(\tau,\t x_2(\tau;{\bf x}), \t x_3(\tau;{\bf x}))$ through $(x_1,x_2,x_3)$, is defined by the following ordinary differential equations:
\be\label{de12} \left\{\begin{array}{ll}
\f{d \t x_2(\tau;{\bf x})}{d\tau}=\t\beta_2(\tau,\t x_2(\tau;{\bf x}),\t x_3(\tau;{\bf x})),\\
\f{d \t x_3(\tau;{\bf x})}{d\tau}=\t\beta_3(\tau,\t x_2(\tau;{\bf x}),\t x_3(\tau;{\bf x})),\\
\t x_2(x_1;{\bf x})=x_2,\\
\t x_2(x_1;{\bf x})=x_3.
\end{array}\right. \ee
Since $(\t\beta_2,\t \beta_3)\in C^{2,\alpha}(\O_e)$, it is easy to prove that $(\t x_2(0;{\bf x}),\t x_3(0;{\bf x}))$ belong to $C^{2,\alpha}(\O_e)$ with respect to ${\bf x}$. Indeed, we have the following estimates:
\be\label{de13}
\|(\t x_2(0;{\bf x}),\t x_3(0;{\bf x}))\|_{C^{2,\alpha}(\O_e)}\leq C_3.
\ee
Then by the characteristic method, it holds that:
\be\label{de14} \left\{\begin{array}{ll}
\beta_2({\bf x})=\epsilon\beta_{2}^{in}(\t x_2(0;{\bf x}),\t x_3(0;{\bf x}))+\int_0^{x_1}\bigg(\f{c^2(\t\rho)}{{\t u_1}^2}(\t\beta_2\p_1\bar s-\p_2\bar s)\bigg)(\tau,\t x_2(\tau;{\bf x}),\t x_3(\tau;{\bf x}))d\tau,\\
\beta_3({\bf x})=\epsilon\beta_{3}^{in}(\t x_2(0;{\bf x}),\t x_3(0;{\bf x}))+\int_0^{x_1}\bigg(\f{c^2(\t\rho)}{{\t u_1}^2}(\t\beta_3\p_1\bar s-\p_3\bar s)\bigg)(\tau,\t x_2(\tau;{\bf x}),\t x_3(\tau;{\bf x}))d\tau.
\end{array}\right. \ee
These enable one to obtain the following estimate:
\be\label{de15}
\|(\beta_2({\bf x}),\beta_3({\bf x}))\|_{C^{2,\alpha}(\O_e)}\leq C_4[\epsilon \|(\beta_{2}^{in},\beta_{3}^{in})\|_{C^{2,\alpha}(\mathrm{T}^2)}+\|\bar s\|_{C^{3,\alpha}(\bar \O)}]\leq C_5(\delta^2+\epsilon).
\ee

Step 3. To resolve $\bar B$. Solving
\be\label{de151} \left\{\begin{array}{ll}
\p_1\bar B+\t\beta_2\p_2\bar B+\t\beta_3\p_3\bar B=-\beta_2\p_2B_0-\beta_3\p_3B_0,\\
 \bar B(0,x_2,x_3)=\epsilon B^{in}(x_2,x_3).
\end{array}\right. \ee
As above, $\bar B$ has the following formulae
$$\bar B=\epsilon B^{in}(\t x_2(0;{\bf x}),\t x_3(0;{\bf x}))+\int_0^{x_1}\bigg(-\beta_2\p_2B_0-\beta_3\p_3B_0\bigg)(\tau,\t x_2(\tau;{\bf x}),\t x_3(\tau;{\bf x}))d\tau.$$ Furthermore, the following estimate holds:
\be\label{de152}
\|\bar B\|_{C^{2,\alpha}(\O_e)}\leq \epsilon \|B^{in}\|_{C^{2,\alpha}(\mathrm{T}^2)}+\t C_5\|(\beta_2({\bf x}),\beta_3({\bf x}))\|_{C^{2,\alpha}(\O_e)}\leq C_6(\delta^2+\epsilon).
\ee

Collecting all the estimates (\ref{de10}), (\ref{de15}) and (\ref{de152}) gives
\be\label{de153}
\|\bar s\|_{C^{3,\alpha}(\O_e)}\leq C_7(\delta^2+\epsilon),\|\beta_2,\beta_3,\bar B\|_{C^{2,\alpha}(\O_e)}\leq C_7(\delta^2+\epsilon).
\ee
Here $C_7=max\{C_2,C_5,C_6\}$, which depends only on the background solution, $\O_e$ and $(\beta_{2}^{in},\beta_{3}^{in},B^{in}, s_e)$.

Choose $\epsilon_1$ small enough, such that if $\epsilon\leq \epsilon_1$, then $C_7^2\epsilon <\f{1}{4}$. Set $\delta=2C_7\epsilon$, then $C_7\delta<\f{1}{2}$ and $C_7(\delta^2+\epsilon)<C_7\epsilon+\f{1}{2}\delta=\delta$. This implies that $\Lambda$ maps $\Xi$ to itself.

It remains to show that the mapping $\Lambda: \Xi\rightarrow \Xi$ is a contraction operator. Suppose $\Lambda:(\t{\bar s}^k,\t \beta_2^k,\t \beta_3^k,\t {\bar B}^k)\mapsto (\bar s^k, \beta_2^k, \beta_3^k,\bar B^k)$ for $k=1,2$. Define the difference $(Y_1,Y_2, Y_3, Y_4)=(\bar s^1-\bar s^2,\beta_2^1-\beta_2^2,\beta_3^1-\beta_3^2, \bar B^1-\bar B^2)$ and $(\t Y_1,\t Y_2,\t Y_3,\t Y_4)=(\t{\bar {s}}^1-\t{\bar s}^2,\t \beta_2^1-\t\beta_2^2,\t\beta_3^1-\t \beta_3^2, \t{\bar B}^1-\t{\bar B}^2)$.

Step 1. Estimate of $Y_1$.

We rewrite the equation satisfied by $\bar s$ as follows:
\be\label{de16}
\sum_{i,j=1}^3\p_i\bigg(a_{ij}(\t {\bar s}, \t \beta_2,\t\beta_3,\t{\bar B})\p_j \bar s\bigg)=\bigg(\sum_{j=1}^3a_{1j}(\t {\bar s}, \t\beta_2,\t\beta_3,\t{\bar B})\p_j\t{\bar s}\bigg)^2-\bigg((\p_2\t\beta_2)^2+(\p_3\t\beta_3)^2+2\p_2\t\beta_3\p_3\t\beta_2\bigg).
\ee

Then $Y_1$ satisfies the following elliptic system:
\be\label{de17} \left\{\begin{array}{ll}
\sum_{i,j=1}^3\p_i\bigg(a_{ij}(\t {\bar s}^1, \t\beta_2^1,\t\beta_3^1,\t{\bar B}^1)\p_j Y_1\bigg)=-\sum_{i=1}^3\p_iG_i+H_1-H_2,\\
-\sum_{j=1}^3\bigg(a_{1j}(\t {\bar s}^1, \t\beta_2^1,\t\beta_3^1,\t{\bar B}^1)\p_j Y_1\bigg)(0,x_2,x_3)=G_1,\\
Y_1(1,x_2,x_3)=0.
\end{array}\right. \ee

Here
\be \left\{\begin{array}{ll}
G_i=\sum_{j=1}^3\bigg([a_{ij}(\t {\bar s}^1, \t\beta_2^1,\t\beta_3^1,\t{\bar B}^1)-a_{ij}(\t{\bar s}^2, \t\beta_2^2,\t\beta_3^2,\t{\bar B}^2)]\p_j \bar s^2\bigg), i=1,2,3;\\
H_1=\bigg(\sum_{i=1}^3 a_{1j}(\t {\bar s}^1, \t\beta_2^1,\t\beta_3^1,\t{\bar B}^1)\p_j\t{\bar s}^1\bigg)^2-\bigg(\sum_{i=1}^3a_{1j}(\t{\bar s}^2, \t\beta_2^2,\t\beta_3^2,\t{\bar B}^2)\p_j \t{\bar s}^2\bigg)^2,\\
H_2=\bigg((\p_2\t\beta_2^1)^2+(\p_3\t\beta_3^1)^2+2\p_2\t\beta_3^1\p_3\t\beta_2^1\bigg)-\bigg((\p_2\t\beta_2^2)^2+(\p_3\t\beta_3^2)^2+2\p_2\t\beta_3^2\p_3\t\beta_2^2\bigg).
\end{array}\right.\ee

Thus the following estimate holds:
\be\begin{aligned}\label{de18}
\|Y_1\|_{C^{2,\alpha}(\O_e)}&\leq C_8\bigg(\sum_{i=1}^3\|G_i\|_{C^{1,\alpha}(\O_e)}+\|(H_1,H_2)\|_{C^{\alpha}(\O_e)}+\|G_1\|_{C^{1,\alpha}(\mathrm{T}^2)}\bigg)\\ &\leq C_9\delta \bigg(\|\t Y_1\|_{C^{2,\alpha}(\O_e)}+\|\t Y_2,\t Y_3,\t Y_4\|_{C^{1,\alpha}(\O_e)}\bigg).
\end{aligned}\ee

Step 2. Estimate of $Y_2,Y_3$.

It follows from (\ref{de11}), $Y_2$ and $Y_3$ satisfy the following system:
\be\label{de19} \left\{\begin{array}{ll}
\p_1 Y_2+\t \beta_2^1\p_2Y_2+\t\beta_3^1\p_3Y_2-\f{c^2(\t\rho^1)}{(\t u_1^1)^2}(\t\beta_2^1\p_1-\p_2)Y_1=K_1,\\
\p_1 Y_3+\t \beta_2^1\p_2Y_3+\t\beta_3^1\p_3Y_3-\f{c^2(\t\rho^1)}{(\t u_1^1)^2}(\t\beta_3^1\p_1-\p_3)Y_1=K_2.\\
Y_2(0,x_2,x_3)=0,\\
Y_3(0,x_2,x_3)=0.
\end{array}\right. \ee
Here

\be \left\{\begin{array}{ll}
K_1=-(\t Y_2\p_2\t\beta_2^2+\t Y_3\p_3\t\beta_2^2)+\bigg(\f{c^2(\t\rho^1)}{(\t u_1^1)^2}\t\beta_2^1-\f{c^2(\t\rho^2)}{(\t u_1^2)^2}\t\beta_2^2\bigg)\p_1\bar s^2-\bigg(\f{c^2(\t\rho^1)}{(\t u_1^1)^2}-\f{c^2(\t\rho^2)}{(\t u_1^2)^2}\bigg)\p_2\bar s^2,\\
K_2=-(\t Y_2\p_2\t\beta_3^2+\t Y_3\p_3\t\beta_3^2)+\bigg(\f{c^2(\t\rho^1)}{(\t u_1^1)^2}\t\beta_3^1-\f{c^2(\t\rho^2)}{(\t u_1^2)^2}\t\beta_3^2\bigg)\p_1\bar s^2-\bigg(\f{c^2(\t\rho^1)}{(\t u_1^1)^2}-\f{c^2(\t\rho^2)}{(\t u_1^2)^2}\bigg)\p_3\bar s^2.
\end{array}\right.\ee

Then it follows directly that:
\be\begin{aligned}\label{de20}
\|(Y_2, Y_3)\|_{C^{1,\alpha}(\O_e)}&\leq C_{10}\bigg(\delta \|(\t Y_1,\t Y_2, \t Y_3, \t Y_4)\|_{C^{1,\alpha}(\O_e)}+\|\nabla Y_1\|_{C^{1,\alpha}(\O_e)}\bigg)\\&\leq C_{11}\delta\bigg(\|\t Y_1\|_{C^{2,\alpha}(\O_e)}+\|\t Y_2,\t Y_3,\t Y_4\|_{C^{1,\alpha}(\O_e)}\bigg).
\end{aligned}\ee

Step 3. Estimate of $Y_4$.

Indeed, $Y_4$ satisfies the following equation:
\be\label{de21} \left\{\begin{array}{ll}
\p_1\bar Y_4+\t\beta_2^1\p_2Y_4+\t\beta_3^1\p_3Y_4=-\t Y_2\p_2\bar B^2-\t Y_3\p_3\bar B^2-Y_2\p_2B_0-Y_3\p_3B_0,\\
 \bar Y_4(0,x_2,x_3)=0.
\end{array}\right. \ee
Then we have the following estimate:
\be\begin{aligned}\label{de22}
\|Y_4\|_{C^{1,\alpha}(\O_e)}&\leq C_{12}\delta \|(\t Y_2, \t Y_3)\|_{C^{1,\alpha}(\O_e)}+C_{12}\|(Y_2, Y_3)\|_{C^{1,\alpha}(\O_e)}\\&\leq C_{13}\delta\bigg(\|\t Y_1\|_{C^{2,\alpha}(\O_e)}+\|\t Y_2,\t Y_3,\t Y_4\|_{C^{1,\alpha}(\O_e)}\bigg).
\end{aligned}\ee

Set $C_{14}=max\{C_9,C_{11},C_{13}\}$, then take $\epsilon_2$ small enough such that $C_{14}\epsilon_2<1$. Now we set $\epsilon_0=min\{\epsilon_1,\epsilon_2\}$. Then if $0<\epsilon<\epsilon_0$, $\Lambda$ maps $\Xi$ to itself and is contraction in low order norm. Hence $\Lambda$ has a unique fixed point, which will be the solution to (\ref{de5}). We have finished our proof.

There are a few remarks in order.

\begin{remark}
{\it We remark that the vorticity of the background solution can be large.}
\end{remark}

\begin{remark}
{\it The above iteration scheme can not be applied to a 3-D nozzle with general sections. In general, the velocity field $\beta_2$ and $\beta_3$ will lose one order derivative when integrating along the particle path. The delicate nonlinear coupling between the hyperbolic modes and elliptic modes $(\beta_2,\beta_3)$ makes it extremely difficult to develop an effective iteration scheme in a general 3-D nozzle. How to effectively explore the good property of the quantity $W=\p_2\beta_3-\p_3\beta_2+\beta_3\p_1\beta_2-\beta_2\p_1\beta_3$ will be investigated in the forthcoming paper.
}
\end{remark}

\section {A new formulation of the three dimensional incompressible Euler equations} \label{IEulerConserved}\par \

\subsection {A new formulation for the incompressible Euler equations}\label{INewconservedquantity}\par \

We apply the same idea to get a new formulation for the incompressible Euler equations. The original incompressible Euler equations take the following form:
\be\label{dc1} \left\{\begin{array}{ll}
 \p_1{u_1}+\p_2{u_2}+\p_3{u_3}=0,\\
 u_1\p_1{u_1}+u_2\p_2{u_1}+u_3\p_3{u_1}+\p_1p=0,\\
 u_1\p_1{u_2}+u_2\p_2{u_2}+u_3\p_3{u_2}+\p_2p=0,\\
 u_1\p_1{u_2}+u_2\p_2{u_2}+u_3\p_3{u_2}+\p_2p=0.
\ea\right. \ee
Define $\beta_i=\f{u_i}{u_1},i=2,3$ and the Bernoulli's function $B=\f{1}{2}(u_1^2+u_2^2+u_3^2)+p$. The inner product of ${\bf u}$ and the last three equations yields the Bernoulli's law:
 \be\label{dc2}
u_1\p_1{B}+u_2\p_2{B}+u_3\p_3{B}=0.
\ee

Multiplying the first equation in (\ref{dc1}) by $\f{-u_{i-1}}{u_1^2}$, dividing the i-th equation in (\ref{dc1}) by $u_1^2$ and adding them together for $i=2,3,4$, we obtain the following new system:
\be\label{dc3} \left\{\begin{array}{ll}
 \p_1{\beta_2}+\beta_3\p_3{\beta_2}-\beta_2\p_3{\beta_3}+\f{1}{u_1^2}\p_2 p=0,\\
 \p_1{\beta_3}+\beta_2\p_2{\beta_3}-\beta_3\p_2{\beta_2}+\f{1}{u_1^2}\p_3 p=0,\\
 \p_2\beta_2+\p_3\beta_3-\f{1}{u_1^2}\p_1 p=0,\\
 \p_1{B}+\beta_2\p_2{B}+\beta_3\p_3{B}=0.
\ea\right. \ee

Using the third equation, we rewrite the above system as follows:
\be\label{dc30} \left\{\begin{array}{ll}
\p_2\beta_2+\p_3\beta_3-\f{1}{u_1^2}\p_1 p=0,\\
\p_1{\beta_2}+\beta_2\p_2{\beta_2}+\beta_3\p_3{\beta_2}-\f{\beta_2}{u_1^2}\p_1 p+\f{1}{u_1^2}\p_2 p=0,\\
\p_1{\beta_3}+\beta_2\p_2{\beta_3}+\beta_3\p_3{\beta_3}-\f{\beta_3}{u_1^2}\p_1 p+\f{1}{u_1^2}\p_3 p=0,\\
\p_1{B}+\beta_2\p_2{B}+\beta_3\p_3{B}=0.
\ea\right. \ee

We also carry out a characteristic analysis to get an insight of the structure
of incompressible Euler equations.

We rewrite (\ref{dc30}) in a matrix form like the following
\be\label{dc4} M_1\p_1{\bf U}+M_2\p_2{\bf U}+M_3\p_3{\bf U}=0. \ee Here ${\bf
U}=(p,\beta_2,\beta_3,B)^T$ and
$$
M_1=\left[\begin{array}{rrrr}
-\f{1}{u_1^2} & 0 & 0 & 0\\
-\f{1}{u_1^2}\beta_2 & 1 & 0 & 0\\
-\f{1}{u_1^2}\beta_3 & 0 & 1 & 0\\
0 & 0 & 0 & 1\\
  \end{array}
\right], M_2=\left[\begin{array}{rrrr}
0 & 1 & 0 & 0\\
\f{1}{u_1^2} & \beta_2 & 0 & 0\\
0 & 0 & \beta_2 & 0\\
0 & 0 & 0 & \beta_2\\
  \end{array}
\right],M_3=\left[\begin{array}{rrrr}
0 & 0 & 1 & 0\\
0 & \beta_3 & 0 & 0\\
\f{1}{u_1^2} & 0 & \beta_3 & 0\\
0 & 0 & 0 & \beta_3\\
  \end{array}
\right].
$$

Denote by $\lambda$  the root of $\det(\lambda
M_1-\xi_2M_2-\xi_3M_3)$. A simple calculation shows that
$\det(\lambda
M_1-\xi_2M_2-\xi_3M_3)=-\f{1}{u_1^2}(\lambda-\beta\cdot\xi)^2[\lambda^2+(\xi_2^2+\xi_3^2)]$,
then we find that it has one real root $\lambda_r=\beta\cdot\xi$
with multiplicity $2$ and two conjugate complex roots
$\lambda_c^{\pm}=\pm i\lambda_I$, where
$\lambda_I=\sqrt{\xi_2^2+\xi_3^2}$.
The corresponding left eigenvectors to $\lambda_r$ are the following
$${\bf l}_r^1=(0,\xi_3+\beta_3(\beta\cdot\xi),-\xi_2-\beta_2(\beta\cdot\xi),0), {\bf l}_r^2=(0,0,0,1).$$
The corresponding left eigenvectors to $\lambda_c^{\pm}$ are
$${\bf l}_c^{\pm}={\bf l}_R\pm i{\bf l}_I=(-\beta\cdot\xi,\xi_2,\xi_3,0)\pm i(\sqrt{\xi_2^2+\xi_3^2},0,0,0).$$
The differential operators corresponding to ${\bf l}_R$ and ${\bf
l}_I$ are
$(-(\beta_2\p_2+\beta_3\p_3),\p_2,\p_3,0)$ and
$(1,0,0,0)$. The action of these two differential operations will
give us an elliptic system. We denote by
$\varpi=\p_2\beta_2+\p_3\beta_3$ and find that $\varpi$ and $p$ satisfy the
following elliptic system:
\be\label{dc40} \left\{\ba{l}
-\f{1}{u_1^2}\p_1p+\varpi=0,\\
\p_1 \varpi+\p_2(\f{1}{u_1^2}\p_2p)+\p_3(\f{1}{u_1^2}\p_3p)-\p_2(\f{\beta_2}{u_1^2}\p_1p)-\p_3(\f{\beta_3}{u_1^2}\p_1p)\\ \quad\quad+(\beta_2\p_2+\beta_3\p_3)(\f{1}{u_1^2}\p_1p)+(\p_2\beta_2)^2+(\p_3\beta_3)^2+2\p_2\beta_3\p_3\beta_2=0.
\ea\right. \ee

\begin{remark}{\it
The steady 3-D incompressible Euler system is always an hyperbolic-elliptic coupled system. Here our point of view is different from previous works \cite{Alber} and \cite{TX}. Indeed, in the works of \cite{Alber} and \cite{TX}, the authors regarded the pressure $p$ as the Lagrange multiplier for the divergent free condition.}
\end{remark}

\subsection {The incompressible Euler equations with a constant Bernoulli's function}\label{INewConservationlaws}  \hspace*{\parindent}

In the following, we assume that the Bernoulli's function is a constant and the speed $|{\bf u}|$ has a positive lower bound $\delta$, then we have $u_1^2=\f{2(B-p)}{1+\beta_2^2+\beta_3^2}$. As in the compressible case, the vorticity field will be parallel to the velocity field and $W=\p_2\beta_3-\p_3\beta_2+\beta_3\p_1\beta_2-\beta_2\p_1\beta_3$
will take the place of $\p_2u_3-\p_3u_2$. We define $G=1+\beta_2^2+\beta_3^2$ and a new function $K(p)$ of $p$ satisfying $K'(p)=\f{1}{2(B-p)}$. and then rewrite the equations
 satisfied by $\beta_2$ and $\beta_3$ as follows:
\be\label{dc5} \left\{\ba{l}
G^{-1}(\p_1\beta_2+\beta_2\p_2\beta_2+\beta_3\p_3\beta_2)-\beta_2\p_1K(p)+\p_2K(p)=0,\\
G^{-1}(\p_1\beta_3+\beta_2\p_2\beta_3+\beta_3\p_3\beta_3)-\beta_3\p_1K(p)+\p_3K(p)=0,\\
\ea\right. \ee

Applying $\beta_3\p_1-\p_3$ and $-\beta_2\p_1+\p_2$ to the above two
equations and adding them together, then one can
obtain the equation for $W$:
\be\label{dc6}{\bf D}W+(\p_2\beta_2+\p_3\beta_3)W-G\p_1K(p)W-G^{-1}W{\bf D}G=0.\ee
Due to the first equation in (\ref{dc3}), $W$ satisfies
the following equation:
\be\label{dc7}{\bf D}W-G^{-1}W{\bf D}G=0.\ee
That is
\be\label{dc8}{\bf D}\bigg(\f{W}{G}\bigg)=0.\ee

It is easy to find $G$ satisfies the following
Riccati-type equation:
\be\label{dc11}{\bf D}G-\f{1}{(B-p)}\p_1p
G^2+\f{1}{(B-p)}G{\bf D}p=0.\ee
Define three new variables
$K_1=G^{-\f{1}{2}},K_i=\f{\beta_i}{G^{\f{1}{2}}}=\beta_iK_1,i=2,3$.
Here $K_2$ and $K_3$ behave as the sine function of ``the flow
angles". Define also two new
functions $I(p)=e^{\int_{\underline{p}}^p\f{1}{2t(B-t)}dt}$ and
$Q(p)=\int_{\underline{p}}^{p}\f{1}{2(B-t)}I(t)^{-2}dt$, where ${\underline{p}}$ is a given positive reference pressure. Since the speed $|{\bf u}|$ has a lower bound  $\delta$, it is easy to verify that the above two integrals are convergent, so $I(p)$ and $Q(p)$ are well-defined. One should note that
$I(p)$ satisfies $\f{I'(p)}{I(p)}=\f{1}{2(B-p)}$.

It follows from (\ref{dc11}) that $K_1$ satisfies:
\be\label{dc12}{\bf D}K_1-\f{1}{2(B-p)}K_1{\bf D}p+\f{1}{2(B-p)}K_1^{-1}\p_1p
=0.\ee
This equation, together with (\ref{dc30}) yields:
\be\label{dc13} \left\{\ba{l}
-\f{1}{2(B-p)}K_1{\bf D}p+\p_1K_1+\p_2K_2+\p_3K_3=0,\\
{\bf D}K_1-\f{1}{2(B-p)}K_1{\bf D}p+\f{1}{2(B-p)}K_1^{-1}\p_1p
=0.\\
{\bf D}K_2-\f{1}{2(B-p)}K_2{\bf D}p+\f{1}{2(B-p)}K_1^{-1}\p_2p=0,\\
{\bf D}K_3-\f{1}{2(B-p)}K_3{\bf D}p+\f{1}{2(B-p)}K_1^{-1}\p_3p=0.
\ea\right. \ee

It follows from the definition of $I(p)$ and $Q(p)$ that

\be\label{dc14} \left\{\ba{l}
-\f{1}{I(p)}K_1{\bf D}I(p)+\p_1K_1+\p_2K_2+\p_3K_3=0,\\
{\bf D}K_1-\f{1}{I(p)}K_1{\bf D}I(p)+\f{1}{I(p)}K_1^{-1}\p_1I(p)
=0.\\
{\bf D}K_2-\f{1}{I(p)}K_2{\bf D}I(p)+\f{1}{I(p)}K_1^{-1}\p_2I(p)=0,\\
{\bf D}K_3-\f{1}{I(p)}K_3{\bf D}I(p)+\f{1}{I(p)}K_1^{-1}\p_3I(p)=0.
\ea\right. \ee

Setting $A_i=\f{K_i}{I(p)}, i=1,2,3$, we
obtain the following conservation laws:
\be\label{dc15} \left\{\ba{l}
\p_1A_1+\p_2A_2+\p_3A_3=0,\\
\p_1( A_1^2)+\p_2( A_1A_2)+\p_3(A_1A_3)+\p_1 Q(p)=0,\\
\p_1( A_1A_2)+\p_2( A_2^2)+\p_3( A_2A_3)+\p_2 Q(p)=0,\\
\p_1( A_1A_3)+\p_2( A_2A_3)+\p_3( A_3^2)+\p_3 Q(p)=0.
\ea\right. \ee

For the incompressible Euler equations, $W$ satisfies the following equation:
\be\label{dc16}(\p_1+\beta_2\p_2+\beta_3\p_3)\bigg(\f{W}{G}\bigg)+\f{1}{2(B-p)^2}[(1,\beta_2,\beta_3)\cdot(\nabla p\times\nabla B)]=0.\ee
And the following conservation law also holds:
\be\label{dc16}
\p_1\bigg(\f{K_1}{I(p,B)}\bigg)+\p_2\bigg(\f{K_2}{I(p,B)}\bigg)+\p_3\bigg(\f{K_3}{I(p,B)}\bigg)=0.
\ee

\begin{remark}{\it
As an application, one can apply the techniques developed in section \ref{Subsoniceuler} to construct a smooth incompressible Euler flow in a rectangular cylinder, which satisfies the given incoming flow angles and the Bernoulli's function at the inlet and the end pressure at the exit. Indeed, our background solution will be denoted by $(p_0,u_0(x_2,x_3),0,0)$. The corresponding Bernoulli's function will be $B_0(x_2,x_3)=\frac{1}{2}u_0^2(x_2,x_3)+p_0$. At the exit of the nozzle, we will prescribe the end pressure:
\begin{equation}\label{dc17} p(1,x_2,x_3)=p_0+\epsilon p_e(x_2,x_3).
\end{equation}
Here it is required that $p_e(x_2,x_3)$ should satisfy the same compatibility conditions as (\ref{de101}) for $s_e(x_2,x_3)$ in the previous section. For the incoming flow angles and the Bernoulli's function and the slip boundary conditions, we should impose the same compatibility conditions as the compressible case. Then we can obtain the following existence and uniqueness results.
}
\end{remark}

\begin{theorem}\label{IETH}
{\it Given $(\beta_{2}^{in},\beta_{3}^{in},B^{in}, p_e)\in C^{3,\alpha}(\mathrm{T}^2)$ satisfying the same compatibility conditions as  (\ref{de100}) and (\ref{de101}) for the compressible case, there exists a positive small number $\epsilon_0$, which depends on the background subsonic state $(p_0,u_0(x_2,x_3),0,0)$ and $(\beta_{2}^{in},\beta_{3}^{in},B^{in}, p_e)$, such that if $0<\epsilon<\epsilon_0$, then there exists a unique smooth subsonic flow $(p,u_1,u_2,u_3)$$\in C^{3,\alpha}(\O_e)\times (C^{2,\alpha}(\O_e))^3 $ to (\ref{dc1}) satisfying boundary conditions (\ref{de1}),(\ref{dc17}) and (\ref{de3}). Moreover, the following estimate holds:
\be\label{de50}
\|(u_1,u_2,u_3)-(u_0(x_2,x_3),0,0)\|_{C^{2,\alpha}(\O_e)}+\|p-p_0\|_{C^{3,\alpha}(\O_e)}\leq C\epsilon.
\ee
Here $C$ is a constant depending on $(p_0,u_0(x_2,x_3),0,0)$ and $(\beta_{2}^{in},\beta_{3}^{in},B^{in}, p_e)$.
}
\end{theorem}

\section {Appendix}\label{appendix}\hspace*{\parindent}

In this appendix, we give the detailed calculations for (\ref{da6})-(\ref{da8}).

Applying $-\beta_3\p_1+\p_3$ and $\beta_2\p_1-\p_2$ to (\ref{da5}) and adding them together, we obtain

\be\begin{aligned}
0&=(\beta_3\p_1-\p_3)(G^{-1}{\bf D}\beta_2)-(\beta_2\p_1-\p_2)(G^{-1}{\bf D}\beta_3)\\&\quad\quad-\bigg[(\beta_3\p_1-\p_3)(\beta_2\p_1-\p_2)K(s)-(\beta_2\p_1-\p_2)(\beta_3\p_1-\p_3)K(s)\bigg]\\&
 =(\beta_3\p_1-\p_3)(G^{-1}{\bf D}\beta_2)-(\beta_2\p_1-\p_2)(G^{-1}{\bf D}\beta_3)-W \p_1 K(s)\\&= J-W\p_1 K(s).
\end{aligned}\ee

While
\be
\begin{aligned}
J=&G^{-1}{\bf D}((\beta_3\p_1-\p_3)\beta_2-(\beta_2\p_1-\p_2)\beta_3)\\&+G^{-1}\sum_{j=1}^3((\beta_3\p_1-\p_3)\beta_j\p_j\beta_2-(\beta_2\p_1-\p_2)\beta_j\p_j\beta_3)\\&
   +G^{-1}\sum_{j=1}^3\beta_j\bigg\{[(\beta_3\p_1-\p_3)\p_j-\p_j(\beta_3\p_1-\p_3)]\beta_2
   \\&\quad\quad\quad\quad-[(\beta_2\p_1-\p_2)\p_j-\p_j(\beta_2\p_1-\p_2)]\beta_3\bigg\}\\&
   +(\beta_3\p_1-\p_3)G^{-1}{\bf D}\beta_2-(\beta_2\p_1-\p_2)G^{-1}{\bf D}\beta_3\\ = &G^{-1}{\bf D}W+J_1+J_2+J_3.
\end{aligned}\ee

Now we compute $J_i, i=1,2,3$ respectively.

\be\begin{aligned}
J_1&=G^{-1}\bigg[W(\p_2\beta_2+\p_3\beta_3)-(\beta_3\p_1-\p_3)\beta_2\p_3\beta_3+(\beta_3\p_1-\p_3)\beta_3\p_3\beta_2
\\&\quad\quad\quad-(\beta_2\p_1-\p_2)\beta_2\p_2\beta_3+(\beta_2\p_1-\p_2)\beta_3\p_2\beta_2\bigg]\\&
= G^{-1}[W(\p_2\beta_2+\p_3\beta_3)+Z].
\end{aligned}\ee

Here $Z=\beta_3\p_1\beta_3\p_3\beta_2-\beta_3\p_1\beta_2\p_3\beta_3+\beta_2\p_1\beta_3\p_2\beta_2-\beta_2\p_1\beta_2\p_2\beta_3$.

\be\begin{aligned}
J_2&=G^{-1}\sum_{j=1}^3\beta_j\bigg\{[(\beta_3\p_1-\p_3)\p_j-\p_j(\beta_3\p_1-\p_3)]\beta_2\\&\quad\quad\quad\quad-[(\beta_2\p_1-\p_2)\p_j-\p_j(\beta_2\p_1-\p_2)]\beta_3\bigg\}\\&
    =G^{-1}Z.
\end{aligned}\ee

\be\begin{aligned}
J_3&=-G^{-2}\bigg[(\beta_3\p_1-\p_3)G{\bf D}\beta_2-(\beta_2\p_1-\p_2)G{\bf D}\beta_3\bigg]\\&=-2G^{-2}\sum_{j=1}^3\bigg[\beta_j(\beta_3\p_1-\p_3)\beta_j{\bf D}\beta_2-\beta_j(\beta_2\p_1-\p_2)\beta_j{\bf D}\beta_3\bigg]\\&
    = -2G^{-2}[W(\beta_2{\bf D}\beta_2+\beta_3{\bf D}\beta_3)+J_{31}].
\end{aligned}\ee
Here
\be\begin{aligned}
J_{31}&=\beta_2(\beta_2\p_1-\p_2)\beta_3{\bf D}\beta_2+\beta_3(\beta_3\p_1-\p_3)\beta_3{\bf D}\beta_2\\&\quad\quad-\beta_3(\beta_3\p_1-\p_3)\beta_2{\bf D}\beta_3-\beta_2(\beta_2\p_1-\p_2)\beta_2{\bf D}\beta_3\\&
       =\beta_2\bigg[(\beta_2\p_1-\p_2)\beta_3{\bf D}\beta_2-(\beta_2\p_1-\p_2)\beta_2{\bf D}\beta_3\bigg]\\&\quad\quad+\beta_3\bigg[(\beta_3\p_1-\p_3)\beta_3{\bf D}\beta_2-(\beta_2\p_1-\p_2)\beta_2{\bf D}\beta_3\bigg]\\&
       =\beta_2\bigg[(\beta_2\p_1-\p_2)\beta_3\p_1\beta_2-(\beta_2\p_1-\p_2)\beta_2\p_1\beta_3\bigg]
       \\&\quad+\beta_3\bigg[(\beta_3\p_1-\p_3)\beta_3\p_1\beta_2-(\beta_3\p_1-\p_3)\beta_2\p_1\beta_3\bigg]\\&
        \quad+\bigg[(\beta_2\p_1-\p_2)\beta_3(\beta_2^2\p_2\beta_2+\beta_2\beta_3\p_3\beta_2)-(\beta_2\p_1-\p_2)\beta_2(\beta_2^2\p_2\beta_3+\beta_2\beta_3\p_3\beta_3)\bigg]\\&
        \quad+\bigg[(\beta_3\p_1-\p_3)\beta_3(\beta_2\beta_3\p_2\beta_2+\beta_3^2\p_3\beta_2)-(\beta_3\p_1-\p_3)\beta_2(\beta_2\beta_3\p_2\beta_3+\beta_3^2\p_3\beta_3)\bigg]\\&  =GZ.
\end{aligned}\ee

Hence we have $J_3=-2G^{-2}[W(\beta_2D\beta_2+\beta_3D\beta_3)]-2G^{-1}Z$. Substitute these calculations into the above formula to get:

\be\begin{aligned}
0&=J-W\p_1 K(s)=J_1+J_2+J_3+G^{-1}DW-W\p_1K(s)\\&
 =G^{-1}[W(\p_2\beta_2+\p_3\beta_3)+Z]+G^{-1}Z-2G^{-2}[W(\beta_2D\beta_2+\beta_3D\beta_3)]\\&\quad\quad-2G^{-1}Z+G^{-1}DW-W\p_1K(s)\\&
 =G^{-1}DW+G^{-1}(\p_2\beta_2+\p_3\beta_3-G\p_1K(s))W-G^{-2}W DG\\&
 =G^{-1}DW-G^{-1}WDs-G^{-2}W DG=\rho D(\f{W}{\rho G}).
\end{aligned}\ee

This implies that $(\p_1+\beta_2\p_2+\beta_3\p_3)(\f{W}{\rho G})=0$.

{\bf Acknowledgements.} This is part of
Ph.D thesis of the author written under the supervision of Professor Zhouping Xin
at the Chinese University of Hong Kong. The author would like to express his gratitude to Prof. Zhouping Xin for his encouragement and wonderful discussion. The author also thanks The Institute of Mathematical Science at The Chinese University of Hong Kong for the support during his Ph.D study. The author would also thank Dr. Wei Yan for the discussion during the author attended the ICCM 2010.

\bibliographystyle{plain}

\end{document}